\documentclass{jams-l}
\usepackage{times}
\usepackage{amsfonts}
\usepackage{amsmath}
\usepackage{amscd}
\usepackage{amssymb}
\usepackage{latexsym}

\theoremstyle{plain}

\swapnumbers
\newtheorem{theorem}[subsection]{Theorem}
\newtheorem{corollary}[subsection]{Corollary}
\newtheorem{lemma}[subsection]{Lemma}
\newtheorem{proposition}[subsection]{Proposition}

\theoremstyle{definition}
\newtheorem{definition}[subsection]{Definition}

\newtheorem{example}[subsection]{Example}

\theoremstyle{remark}

\newtheorem{remark}[subsection]{Remark}

\swapnumbers

\newcommand\preprintnote {preprint on \myhomepage}
\newcommand\myhomepage{http://www.math.ohio-state.edu/\-\~{}schoutens}

\newcommand\OSU{\address{Department of Mathematics\\
100 Math Tower\\
Ohio State University\\
Columbus, OH 43210 (USA)}
\email{schoutens@math.ohio-state.edu}}

\newcommand{\emptyprop}{q}

\newcommand \complet[1]{\widehat {#1}} 
\newcommand \exactseq [5]{0\to{#1}\:\map{#2}\:{#3}\:\map{#4}\:{#5}\to0}

\newcommand \ext[4]{\operatorname{Ext}_{#1}^{#2}(#3,#4)}
 
\newcommand \id{\mathfrak a}

\newcommand \iso{\cong}

\newcommand \map[1]{{\newcommand{\tmpprop}{#1q}  \if\tmpprop\emptyprop \to\else \xrightarrow{{\phantom{i}{#1}\phantom{i}}}\fi}} 
\newcommand \maxim{\mathfrak m}
\newcommand \nat{\mathbb N}
\newcommand \norm[1]{\left|#1\right|}

\newcommand \pol[2]{#1[#2]}
\newcommand \pow[2]{#1[[#2]]}

\newcommand \pr{\mathfrak p}

\newcommand \range [2]{#1,\dots,#2}

\newcommand \rij[2]{(#1_1,\dots,#1_{#2})}

\newcommand \tensor{\otimes}
\newcommand \tor[4]{\operatorname{Tor}^{#1}_{#2}(#3,#4)}
\newcommand \op\operatorname

\newcommand{\commdiagram}[9][]{%
\begin{equation}
{\newcommand{\tmpprop}{#1q} 
\if\tmpprop\emptyprop \relax\else \label{#1}\fi}
\begin{aligned}%
\mbox{
\begin{picture}(130,90)%
\put(120,70){\vector( 0,-1){50}}%
\put(10,80){\vector( 1, 0){100}}%
\put(0,70){\vector( 0,-1){50}}%
\put(10,10){\vector( 1, 0){100}}%
\put(115,80){\makebox(0,0)[l]{$#4$}}%
\put(5,80){\makebox(0,0)[r]{$#2$}}%
\put(115,10){\makebox(0,0)[l]{$#9$}}%
\put(5,10){\makebox(0,0)[r]{$#7$}}%
\put(-3,50){\makebox(0,0)[r]{$#5$}}
\put(123,50){\makebox(0,0)[l]{$#6$}}
\put(60,3){\makebox(0,0)[c]{$#8$}}
\put(60,88){\makebox(0,0)[c]{$#3$}}
\end{picture}}
\end{aligned}
\end{equation}}

\newcommand\commtrianglefront[7][]{%
\begin{equation}
{\newcommand{\tmpprop}{#1q} 
\if\tmpprop\emptyprop \relax\else \label{#1}\fi}
\begin{aligned}%
\mbox{
\begin{picture}(120,80)%
\put(55,70){\vector(-1,-2){30}}
\put(65,70){\vector(1,-2){30}}
\put(30,5){\vector(1,0){60}}
\put(60,75){\makebox(0,0)[c]{$#2$}}
\put(25,5){\makebox(0,0)[r]{$#4$}}
\put(95,5){\makebox(0,0)[l]{$#6$}}
\put(60,0){\makebox(0,0)[c]{$#5$}}
\put(37,43){\makebox(0,0)[r]{$#3$}}
\put(83,43){\makebox(0,0)[l]{$#7$}}
\end{picture}}
\end{aligned}
\end{equation}}

\newcommand\commtriangleback[7][]{%
\begin{equation}
{\newcommand{\tmpprop}{#1q}
\if\tmpprop\emptyprop \relax\else \label{#1}\fi}
\begin{aligned}%
\mbox{
\begin{picture}(120,80)%
\put(55,70){\vector(-1,-2){30}}
\put(65,70){\vector(1,-2){30}}
\put(30,5){\vector(1,0){60}}
\put(60,75){\makebox(0,0)[c]{$#2$}}
\put(25,5){\makebox(0,0)[r]{$#6$}}
\put(95,5){\makebox(0,0)[l]{$#4$}}
\put(60,0){\makebox(0,0)[c]{$#5$}}
\put(37,43){\makebox(0,0)[r]{$#7$}}
\put(83,43){\makebox(0,0)[l]{$#3$}}
\end{picture}}
\end{aligned}
\end{equation}}

\newcommand{\name}[1]{{\sc#1}}
\newcommand \acf{algebraically closed field}

\newcommand \ch{characteristic}
\newcommand \homo{homomorphism}
\newcommand \CM{Coh\-en-Mac\-au\-lay}
\renewcommand\iff{if, and only if,}
\newcommand \DVR{discrete valuation ring}

\DeclareMathOperator*{\UP}{ulim}
\newcommand \up[1]{\UP_{#1\to\infty}}

\newcommand \ul[1]{\seq{#1}\infty}
\newcommand \seq[2]{#1\mathstrut_{#2}}
\newcommand \sr{approximation}
\newcommand \SR{Approximation}
\newcommand \uleq[1]{\ulseq{#1}\infty{eq}}
\newcommand \ulseq[3]{#1\mathstrut^{\text{#3}}_{#2}}
\newcommand \ulmix[1]{\ulseq{#1}\infty{mix}}
\newcommand \los{\L os' Theorem}

\newcommand \fitt[1]{\Re_{#1}}
\newcommand \br{\mathfrak O}
\newcommand \affine[1]{\underline{\textsf{Aff}}(#1)}

\newcommand  \q{pseudo}
\newcommand  \Q{Pseudo}

\newcommand  \qdim{\q-dimension}
\newcommand  \Qdim{\Q-dimension}

\newcommand  \qCM{\q-\CM}
\newcommand  \qreg{\q-regular}

\newcommand  \qht{\q-height}
\newcommand  \Qht{\Q-height}

\newcommand  \serene{isodimensional}
\newcommand  \blim{restricted ultraproduct}
\newcommand  \Blim{Restricted Ultraproduct}
\newcommand  \sepquot{separated quotient}
\newcommand  \BCM[1]{\mathcal B(#1)}

\newcommand \bound[2]{\textsf{#1}(#2)}

\newcommand \zet{\mathbb Z}
\renewcommand{\name}[1]{#1}

\title {Asymptotic Homological Conjectures in mixed characteristic}
\author{Hans Schoutens
}
\thanks{Partially supported by a  grant from the National Science Foundation.}
\date{30.03.2003}
\OSU
\urladdr{\myhomepage} 

\begin{document}

\begin{abstract}   
In this paper, various Homological Conjectures are studied for local rings which are  locally finitely generated over a \DVR\ $V$ of mixed \ch. Typically, we can only conclude that a particular Conjecture holds for such a ring provided the residual \ch\ of $V$ is sufficiently large in terms of the complexity of the data, where the complexity is primarily given in terms of the degrees of the polynomials over $V$ that define the data, but possibly  also by some additional invariants such as (homological) multiplicity. Thus asymptotic versions of the Improved New Intersection Theorem, the Monomial Conjecture, the Direct Summand Conjecture, the Hochster-Roberts Theorem and the Vanishing of Maps of Tors Conjecture  are given. 

That the results only hold  asymptotically, is due to the fact that non-standard arguments are used, relying on the Ax-Kochen-Ershov Principle, to infer their validity from their positive \ch\ counterparts. A key role in this transfer is played by the Hochster-Huneke canonical construction of big \CM\ algebras in positive \ch\ via absolute integral closures. 
\end{abstract}

\maketitle

\section{Introduction}

In the last three decades,  all the so-called Homological Conjectures have been settled completely for Noetherian local rings containing a field by work of \name{Szpiro-Peskine, Hochster-Roberts, Hochster, Evans-Griffith}, et.~al. (some of the main papers are \cite{EG,HoHT,HoDS,HR,PS72}). More recently, \name{Hochster-Huneke} have given more simplified proofs of most of these results by means of their tight closure theory, including their  canonical construction of big \CM\ algebras in positive \ch\ (see \cite{HHbigCM,HH,HHZero,HuTC}; for further discussion and proofs, see \cite[\S9]{BH} or \cite{Str90}).

In sharp contrast is the development in mixed \ch, where only sporadic results (often in low dimensions) are known, apart from a single break-through \cite{Rob87} by \name{Roberts}, settling the New Intersection Theorem for all Noetherian local rings. Some attempts  have been made by \name{Hochster}, either by finding a suitable substitute for tight closure in mixed \ch\  \cite{HoSol},  or by finding big \CM\ modules in mixed \ch\ \cite{HoWitt}.  These approaches have yet to bear fruit and the best result to date in this direction is the existence of  big \CM\ algebras in dimension three \cite{HoBCM3}, which in turn relies on the positive solution of the  Direct Summand Conjecture in dimension three by Heitmann \cite{HeitDS}.

In this paper, we will follow the big \CM\ algebra approach, but instead of trying to work with rings of Witt vectors, we will use the Ax-Kochen-Ershov Principle \cite{AK,Ers65,Ers66}, linking complete \DVR{s} in mixed \ch\ with complete \DVR{s} in positive \ch\ via an equi\ch\ zero (non-discrete) valuation ring (see Theorem~\ref{T:AKE} below). This intermediate valuation ring is obtained by a construction which originates from logic, but is quite algebraic in nature, to wit, the \emph{ultraproduct} construction. Roughly speaking, this construction associates to an infinite collection of rings $\seq Cw$ their ultraproduct $\ul C$, which should be thought of as a kind of ``limit'' or ``average''  (realized as a certain homomorphic image of the product). An ultraproduct inherits many of the algebraic properties of its components. The correct formulation of this transfer principle is \los, which makes precise when a property carries over (namely, when it is first order definable in some suitable language). Properties that carry over are those of being \emph{a domain, a field, a valuation ring, local, Henselian}; among the properties that do not carry over is \emph{Noetherianity}, so that almost no ultraproduct is Noetherian (except an ultraproduct of fields or of Artinian rings of bounded length). This powerful tool is used in \cite{Schm87,SvdD,SchBC,SchBArt}, to obtain uniform bounds in polynomial rings over fields; in \cite{SchBC,SchBounds,SchSymPow,SchBS}, to transfer properties from positive to zero \ch; and in \cite{SchBCM,SchHR,SchNSTC}, to give an alternative treatment of tight closure theory in equi\ch\ zero. The key fact in the first set of papers is a certain flatness result about ultraproducts (see Theorem~\ref{T:SvdD} below for a precise formulation), and in the two last sets, the so-called \emph{Lefschetz Principle for \acf{s}} (the ultraproduct of the algebraic closures of the $p$-element fields $\mathbb F_p$ is isomorphic to $\mathbb C$). 

The Ax-Kochen-Ershov Principle is a kind of \emph{Lefschetz Principle for Henselian valued fields}, and its most concrete form states that the ultraproduct of all $\pow{\mathbb F_p}t$, with $t$ a single variable,  is isomorphic  to the ultraproduct of all rings of $p$-adic integers $\zet_p$. We will identify both ultraproducts and denote the resulting ring by $\br$. It follows that $\br$ is an equi\ch\ zero Henselian valuation ring with principal maximal ideal, whose  separated quotient (=the reduction modulo the intersection of all powers of the maximal ideal) is an equi\ch\ zero excellent  complete \DVR.

To explain  the underlying idea in this paper, for a (not necessarily Noetherian) local ring $(Z,\pr)$, let us denote by $\affine Z$  the category of \emph{local $Z$-affine algebras}, that is to say, $Z$-algebras of the form $(\pol ZX/I)_\maxim$, with $X$ a finite tuple of variables, $I$ a finitely generated ideal and $\maxim$ a prime ideal containing $\pr$ and $I$. The objective is to transfer algebraic properties (such as the homological Conjectures)  from the positive \ch\ categories $\affine{\pow{\mathbb F_p}t}$ to the mixed \ch\ categories $\affine{\zet_p}$. This will be achieved through the intermediate equi\ch\ zero category $\affine\br$. As this latter category consists mainly of non-Noetherian rings, we will have to find analogues in this setting of many familiar notions from commutative algebra, such as dimension, depth, \CM{ness} or regularity (see \S\S\ref{s:qdim} and \ref{s:qsing}). 

The following  example is paradigmatic: let $X$ be a finite tuple of variables and let $\uleq A$ be the ultraproduct of all $\pol{\pow{\mathbb F_p}t}X$, and $\ulmix A$, the ultraproduct of all $\pol{\zet_p}X$. Note that both rings contain $\br$, and in fact, contain $\pol\br X$. The key algebraic fact, which follows from a result on effective bounds by \name{Aschenbrenner}, is that both inclusions are flat. Suppose we have in each $\pol{\pow{\mathbb F_p}t}X$ a polynomial $\seq fp$, and let $\ul f$ be their ultraproduct. A priori, we can only say  that $\ul f\in\uleq A$. However, if all $\seq fp$ have degree $d$, for some $d$ independent from $p$, then $\ul f$ itself is a polynomial over $\br$ of degree $d$ (since an ultraproduct commutes with finite sums by \los). Hence, as $\ul f$ lies in $\pol \br X$, we can also view it as an element in $\ulmix A$. Therefore, there are polynomials $\seq{\tilde f}p\in\pol{\zet_p}X$ whose ultraproduct is equal to $\ul f$. The choice of the $\seq{\tilde f}p$ is not unique, but any two choices will be equal for almost $p$, by \los. In conclusion, to a collection of polynomials defined over the various $\pow{\mathbb F_p}t$, of uniformly bounded degree,  we can associate, albeit not uniquely, a collection of polynomials defined over the various $\zet_p$ (of uniformly bounded degree), and of course, this also works the other way. Instead of doing this for just one polynomial in each component, we can now do this for a finite tuple of polynomials of fixed length. If at the same time, we can maintain certain algebraic relations among them (characterizing one of the properties we seek to transfer), we will have achieved our goal.  

 Unfortunately, it is the nature of an ultraproduct that it only captures the ``average'' property of its components. In the present context, this means that the desired property does not necessarily hold in all $\pol{\zet_p}X$, but only in almost all. In conclusion, we cannot hope for a full solution of the Homological Conjectures by this method, but only an \emph{asymptotic} solution. In view of the above, the following definition is natural.

\subsubsection*{Complexity.}
Let $Z\to C$ be a local \homo\ between Noetherian local rings. We say that $C$ has \emph{$Z$-complexity} at most $c$, if we can write $C$ as $(\pol ZX/I)_\maxim$ (so that $C$ is a member of $\affine Z$), with $X$ a tuple of at most $c$ variables and $\maxim$ a prime ideal containing $I$ and    the maximal ideal of $Z$, such that $I$ and $\maxim$ are generated by polynomials of degree at most $c$.  

An element $r\in C$ is said to have \emph{$Z$-complexity} at most $c$, if $C$ has $Z$-complexity at most $c$ and if there exist polynomials $P,Q\in\pol ZX$ of degree at most $c$ with $Q\notin \maxim$, such that $r=P/Q$ in $C$. Sometimes we might say that  a tuple or a matrix has $Z$-complexity at most $c$, to indicate that each entry has $Z$-complexity at most $c$ and the dimensions are also bounded by $c$.

An ideal $J$ of $C$ has \emph{$Z$-complexity} at most $c$, if $C$ has $Z$-complexity at most $c$ and $J$  is generated by elements of $Z$-complexity at most $c$. A \homo\ $C\to D$ is said to have \emph{$Z$-complexity} at most $c$, if both $C$ and $D$ have $Z$-complexity at most $c$ and the \homo\ is given by sending $X_i$ to an element in $S$ of $Z$-complexity at most $c$ (where $C=(\pol ZX/I)_\maxim$). 

\subsubsection*{Asymptotic Properties.}
Let $\mathcal P$ be a property of Noetherian local rings (possibly involving some additional data). We will use the phrase  \emph{$\mathcal P$  holds asymptotically in mixed \ch}, to express that for each $c$, we can find  a bound $c'$, such that if $V$ is a complete \DVR\ of mixed \ch\ and $C$ a $V$-algebra of $V$-complexity at most $c$ (and a similar bound on the data), then property $\mathcal P$ holds for $C$, provided the \ch\ of the residue field of $V$ is at least $c'$. Sometimes, we have to control some additional invariants in terms of the bound $c$. In this paper, we will prove that  in this sense, many Homological Conjectures hold asymptotically in mixed \ch.

\subsubsection*{A Final Note.}
Its asymptotic nature is the main impediment of the present method to carry  out  \name{Hochster}'s program of obtaining tight closure and   big \CM\ algebras in mixed \ch. For instance, despite the fact that we  are able to define an analogue of a balanced big \CM\ for $\br$-affine domains, this object cannot be realized as an ultraproduct of $\zet_p$-algebras, so that there is no candidate so far for a big \CM\ in mixed \ch. Although I will not pursue this line of thought in this paper, one could also define some non-standard closure operation on ideals in $\br$-affine algebras, but again, such an operation will only partially  descend to any component.

\subsubsection*{Notation.}
A \emph{tuple} $\mathbf x$ over a ring $Z$ is always understood to be finite. Its length is denoted by $\norm{\mathbf x}$ and the ideal it generates is denoted $\mathbf xZ$. When we say that $(Z,\pr)$ is \emph{local}, we mean that $\pr$ is its maximal ideal, but we do not imply that $Z$ has to be Noetherian.

For a survey of the results and methods in this paper, see also \cite{SchMixBCMCR}.

\section{Ultraproducts}

In this preliminary section, I state some generalities about ultraproducts and then briefly review the situation in equi\ch\ zero and the Ax-Kochen-Ershov Principle. The next section lays out the essential tools for conducting the transfer discussed in the introduction, to wit, \sr{s}, \blim{s} and non-standard hulls, whose properties are then studied in \S\S\ref{s:qdim} and \ref{s:qsing}. The subsequent sections contain proofs of various asymptotic results, using these tools.

Whenever we have an infinite index set $W$, we will equip it with some (unnamed) non-principal ultrafilter; ultraproducts will always be taken with respect to this ultrafilter and we will write $\up w \seq Ow$ or simply $\ul O$ for the ultraproduct of objects $\seq Ow$ (this will apply to rings, ideals and elements alike).  A first introduction to ultraproducts, including \los, sufficient to understand the present paper, can be found in \cite[\S2]{SchNSTC}; for a more detailed treatment, see \cite{Hod}.  \los\ states essentially that if a fixed algebraic relation holds among finitely many elements $\seq{f_1}w,\dots,\seq{f_s}w$ in each ring $\seq Cw$, then the same relation holds among their ultraproducts $\ul{f_1},\dots,\ul{f_s}$ in the ultraproduct $\ul C$, and conversely, if such a relation holds in $\ul C$, then it holds in almost all $\seq Cw$. Here \emph{almost all} means ``for all $w$ in a subset  of the index set which belongs to the ultrafilter'' (the idea is that sets belonging to the ultrafilter are \emph{large}, whereas the remaining sets  are \emph{small}).

An immediate, but important application of \los\ is that the ultraproduct of \acf{s} of different prime \ch{s} is an (uncountable) \acf\ of \ch\ zero, and each uncountable \acf\ of \ch\ zero, including $\mathbb C$, can be realized thus. This simple observation, in combination with work of \name{van den Dries} on non-standard polynomials (see below), was exploited in \cite{SchNSTC} to define an alternative version of tight closure for $\mathbb C$-affine algebras, called \emph{non-standard tight closure}.   The ensuing notions of F-regularity and F-rationality have been proven to be more versatile \cite{SchBCM,SchRatSing}, than those defined by \name{Hochster-Huneke} in \cite{HHZero}.

Let me briefly recall the  results in \cite{SvdD,vdD79} on non-standard polynomials mentioned above. Let $\seq Kw$ be fields (of arbitrary \ch) with ultraproduct $\ul K$ (which is again a field by \los). Let $X$ be a fixed finite tuple of variables and set $A:=\pol{\ul K}X$ and $\seq Aw:=\pol{\seq Kw}X$. Let $\ul A$ be the ultraproduct of the $\seq Aw$. As in the example discussed in the introduction, we have a canonical embedding of $A$ inside $\ul A$. In  fact, the following easy observation, valid over arbitrary rings, describes completely the elements in $\ul A$ that lie in $A$ (the proof is straightforward and left to the reader). 

\begin{lemma}\label{L:deg}
Let  $X$ be a finite tuple of variables. Let $\seq Cw$ be rings and let $\ul C$ be their ultraproduct. If $\seq fw$ is a polynomial in $\pol{\seq Cw}X$ of degree at most $c$, for each $w$ and for some $c$ independent from $w$, then their ultraproduct in $\up w\pol{\seq Cw}X$ belongs already to the subring $\pol{\ul C}X$, and conversely, every element in $\pol{\ul C}X$ is obtained in this way.
\end{lemma}

This result also motivates the notion of \emph{complexity} from the introduction. Returning to the \name{Schmidt-van den Dries} results,  the following two properties of the embedding $A\subset \ul A$ do not only imply the uniform bounds from \cite{SvdD,SchBC}, but play also an important theoretical role in the development of non-standard tight closure \cite{SchNSTC}.

\begin{theorem}[\name{Schmidt-van den Dries}]\label{T:SvdD}
The embedding $A\subset \ul A$ is faithfully flat and every prime ideal in $A$ extends to a prime ideal in $\ul A$.
\end{theorem}

To carry out the present program, we have to replace the base fields $\seq Kw$ by complete \DVR{s} $\seq\br w$. Unfortunately, we now have to face the following complications. Firstly, the ultraproduct $\ul\br $ of the $\seq\br w$ is no longer Noetherian, and so in particular the corresponding $A:=\pol{\ul\br} X$ is non-Noetherian. Moreover, the embedding $A\subset \ul A$, where $\ul A$ is now the ultraproduct of the $\seq Aw:=\pol{\seq\br w}X$, although flat (see Theorem~\ref{T:ff} below), is no longer faithfully flat (this is related to Dedekind's problem; see \cite{Asch} or \cite{SchUBS} for details). Moreover, not every prime ideal extends to a prime ideal, and in order to preserve this, we will have to work locally (see Remark~\ref{R:extprim}). 

To address the first of these problems, we will realize $\ul\br$ in two different ways, as an ultraproduct of complete \DVR{s} in positive \ch\ and as an ultraproduct of complete \DVR{s} in mixed \ch, and then pass from one set to the other via $\ul\br$, as explained in the introduction. This is the celebrated Ax-Kochen-Ershov Principle \cite{AK,Ers65,Ers66}, and I will discuss this now. For each $w$, let $\ulseq\br p{mix}$ be a complete \DVR\ of mixed \ch\ with residue field $\seq\kappa p$ of \ch\ $p$. To each $\ulseq\br p{mix}$, we associate a corresponding equi\ch\ complete \DVR\ with the same residue field, by letting 
	\begin{equation}\label{eq:eqDVR}
	\ulseq\br p{eq}:= \pow{\seq\kappa p}t
	\end{equation}
where $t$ is a single variable.

\begin{theorem}[\name{Ax-Kochen-Ershov}]\label{T:AKE}
The ultraproduct of the $\ulseq\br p{eq}$ is isomorphic (as a local ring) with the ultraproduct of the $\ulseq\br p{mix}$.
\end{theorem}

\begin{remark}
As stated, we need to assume the continuum hypothesis. Otherwise, by the Keisler-Shelah Theorem \cite[Theorem 9.5.7]{Hod}, one might need to take further ultrapowers, that is to say, take a larger index set endowed with a (non-$\omega$-complete) non-principal ultrafilter. In order to not complicate the exposition, I will nonetheless make the  set-theoretic assumption, so that our index set can always be taken to be  the set of prime numbers.
\end{remark}

To conclude this section, I state a variant of Prime Avoidance which will be used in the form discussed in the remark following it.

\begin{proposition}[Controlled Ideal Avoidance]
Let $Z$ be a local ring with infinite residue field $\kappa$ and let $C$ be an arbitrary $Z$-algebra. Let $I:=\rij fnC$ be a finitely generated ideal in $C$ and let $\id_1,\dots,\id_t$ be arbitrary ideals in $C$ not containing $I$. If $W$ denotes the $Z$-submodule of $C$ generated by the $f_i$, then we can find $f\in W$ not contained in any of the $\id_j$.
\end{proposition}
\begin{proof}
We induct on the number $t$ of ideals to be avoided. If $t=1$, then some $f_i\notin\id_1$, since $I\not\subset\id_1$. Hence assume $t>1$. By induction, we can find elements $g_i\in W$, for $i=1,2$, which lie outside any $\id_j$ for $j\neq i$. If either $g_1\notin\id_1$ or $g_2\notin\id_2$ we are done, so assume $g_i\in\id_i$. Therefore, every element of the form $g_1+zg_2$ with $z$ a unit in $Z$ does not lie in $\id_1$ nor in $\id_2$. Since $\kappa$ is infinite, we can find $t-1$ units $z_1,z_2,\dots,z_{t-1}$ in $Z$ whose residues in $\kappa$ are all distinct. I claim that at least one of the $g_1+z_ig_2$ lies outside all $\id_j$, so that we found our desired element in $W$. Indeed, if not, then each $g_1+z_ig_2$ lies in one of the $t-2$ ideals $\id_3,\dots,\id_t$, by our previous remark. By the Pigeon Hole Principle, for some $j$ and some $l\neq k$, we have that $g_1+z_kg_2$ and $g_1+z_lg_2$ lie both in $\id_j$. Hence so does their difference $(z_k-z_l)g_2$. However, $z_k-z_l$ is a unit in $Z$, by choice of the $z_i$,  so that $g_2\in\id_j$, contradiction.
\end{proof}

\begin{remark}\label{R:pa}
In particular, if $Z$ is Noetherian and both $C$ and $I$ have $Z$-complexity at most $c$, then we can find an element $f\in I$ of $Z$-complexity at most $c$, outside any finite set of ideals not containing $I$. Indeed, every element in the module $W$ has $Z$-complexity at most $c$.
\end{remark}

\section{\SR{s}, \Blim{s} and Non-standard hulls}

In this section, some general results on ultraproducts of finitely generated algebras over \DVR{s} will be derived.  We start with introducing some general terminology, over arbitrary Noetherian local rings, but once we start proving some non-trivial properties in the next sections, we will specialize to the case that the base rings are \DVR{s}.  For some results in the general case, we refer to \cite{SchBArt,SchUBS}.

For each $w$, we fix a Noetherian local ring $\seq \br w$ and let $\br$ be its ultraproduct. If the $\seq\pr w$ are the maximal ideals of the $\seq \br w$, then their ultraproduct $\pr$ is the maximal ideal of $\br$. We will write  $\varpi$ for the ideal of \emph{infinitesimals} of $\br$, that is to say, the intersection  of all the powers $\pr^k$ (note that in general $\varpi\neq (0)$ and therefore, $\br$ is in particular non-Noetherian).

By saturatedness of ultraproducts, $\br$ is quasi-complete in its $\pr$-adic topology in the sense that any Cauchy sequence has a (non-unique) limit. Hence the completion of $\br$ is $\br/\varpi$ (see also Lemma~\ref{L:comp} below). Moreover, we will assume that all $\seq \br w$ have embedding dimension at most $\epsilon$. Hence so do  $\br$ and   $\br/\varpi$. Since a complete local ring with finitely generated maximal ideal is Noetherian (\cite[Theorem 29.4]{Mats}), we showed that   $\br/\varpi$ is a Noetherian complete local ring. For more details in the case of interest to us, where each $\seq \br w$ is a \DVR\ or a field, see  \cite{BDDL}.

We furthermore fix throughout a tuple of variables $X=\rij Xn$ and, we let $A:=\pol \br X$ and $\seq Aw:=\pol{\seq \br w}X$.

\begin{definition}
The \emph{non-standard hull} of $A$ is by definition the ultraproduct of the $\seq Aw$ and is denoted $\ul A$.
\end{definition}

By \los, we have an inclusion $\br\subset \ul A$. Let us continue to write $X_i$ for the ultraproduct in $\ul A$ of the constant sequence $X_i\in\seq Aw$. By \los, the $X_i$ are algebraically independent over $\br$. In other words,  $A$ is a subring of $\ul A$. In the next section, we will prove the key algebraic property of the extension $A\subset \ul A$ when the base rings $\seq\br w$ are \DVR{s}, to wit, its flatness. We start with extending the notions of non-standard hull  and \sr\ from \cite{SchNSTC}, to arbitrary local $\br$-affine algebras (recall that a \emph{local $\br$-affine algebra} is a localization of a finitely presented $\br$-algebra at a prime ideal containing $\pr$).

\subsubsection*{\SR{s} and non-standard hulls.}
An \emph\sr\ of a polynomial $f\in A$ is a sequence of polynomials $\seq fw\in\seq Aw$, such that their ultraproduct is equal to $f$, viewed as an element in $\ul A$. Note that according to Lemma~\ref{L:deg}, we can always find such an \sr. Moreover, any two \sr{s} are equal for almost all $w$, by \los. Similarly, an \emph{\sr} of a   finitely generated ideal $I:=\mathbf fA$ with $\mathbf f$ a finite tuple, is a sequence of ideals $\seq Iw:=\seq{\mathbf f}w\seq Aw$, where $\seq{\mathbf f}w$ is an \sr\ of $\mathbf f$ (meaning that each entry in $\seq{\mathbf f}w$ is an \sr\ of the corresponding entry in $\mathbf f$). \los\ gives once more that any two \sr{s} are almost all equal. Moreover, if $\seq Iw$ is some \sr\ of $I$ then
	\begin{equation}\label{eq:ulid}
	\up w\seq Iw=I\ul A.
	\end{equation}

Assume now that $C$ is a finitely presented $\br$-algebra, say $C=A/I$ with $I$ a finitely generated ideal. We define an \emph\sr\ of $C$ to be the sequence of finitely generated $\seq \br w$-algebras $\seq Cw:=\seq Aw/\seq Iw$, where $\seq Iw$ is some \sr\ of $I$. We define the \emph{non-standard hull} of $C$ to be the ultraproduct of the $\seq Cw$ and denote it $\ul C$. It is not hard to show that $\ul C$ is uniquely defined up to $C$-algebra isomorphism (for more details see \cite{SchNSTC} or \cite{SchBArt}). From \eqref{eq:ulid}, it follows that $\ul C=\ul A/I\ul A$. In particular, there is a canonical \homo\ $C\to \ul C$ obtained from the base change $A\to \ul A$. 

\subsubsection*{Some Caveats.}
When $I$ is not finitely generated, $I\ul A$ might not be realizable as an ultraproduct of ideals, and consequently, has no \sr. Although one can find special cases of infinitely generated ideals admitting \sr{s}, we will never have to do this in the present paper. Similarly, we only define \sr{s} for finitely presented algebras.

Although $A\to \ul A$ is injective, this is not necessarily the case for $C\to\ul C$, if the $\seq \br w$ are not fields. For instance, if $W$ is the set of prime numbers, $\seq \br p:=\zet_p$ for each $p\in W$ and $I=(1-\pi X,\gamma)A$ where $\pi=\up pp$ and $\gamma=\up pp^p$, then $I\neq (1)$ but $I\ul A=(1)$. 

However, when the $\seq \br w$ are \DVR{s}, we will see shortly, that this phenomenon disappears if we localize at prime ideals containing $\pr$.  Next we define a process which is converse to taking \sr{s}.

\subsubsection*{\Blim{s}.}
Fix some $c$.  For each $w$, let $\seq Iw$ be an ideal in $\seq Aw$ of $\seq \br w$-complexity at most $c$. In other words, we can write $\seq Iw=\seq{\mathbf f}w\seq Aw$, for some tuple $\seq{\mathbf f}w$ of $\seq \br w$-complexity at most $c$. Let $\mathbf f$ be the ultraproduct of these tuples. By Lemma~\ref{L:deg}, the tuple $\mathbf f$ is already defined over $A$. We call $I:=\mathbf fA$ the \emph\blim\ of the $\seq Iw$. It follows that the $\seq Iw$ are an \sr\ of $I$ and that $I\ul A$ is the ultraproduct of the $\seq Iw$. 

With $\seq Cw:=\seq Aw/\seq Iw$   and  $C:=A/I$, we call $C$ the \emph\blim\ of the $\seq Cw$. The $\seq Cw$ are an \sr\ of $C$ and their ultraproduct $\ul C$ is the non-standard hull of $C$. We can now extend the previous definition to the image in $\seq Cw$ of an element $\seq cw\in\seq Aw$ (respectively, to the extension $\seq Jw\seq Cw$ of a finitely generated ideal $\seq Jw\subset \seq Aw$) of $\seq \br w$-complexity at most $c$ and define similarly their \emph\blim\ $c\in C$ and $JC$ as the image in $C$ of the respective \blim\ of the $\seq cw$ and the $\seq Jw$.

\subsubsection*{Functoriality.}
We have a commutative diagram
	\commdiagram C\varphi D {} {} {\ul C} {\ul\varphi} {\ul D}
where $C\to D$ is an $\br $-algebra \homo\ of finite type between finitely presented $\br $-algebras and $\ul C\to \ul D$ is its base change over the non-standard hulls $\ul C:=\ul A\tensor C$ and $\ul D:=\ul A\tensor D$, or alternatively, where $C\to D$ (respectively, $\ul C\to\ul D$) are the \blim\ (respectively, ultraproduct) of $\seq \br w$-algebras \homo{s} $\seq Cw\to\seq Dw$ of  $\seq \br w$-complexity at most $c$, for some $c$ independent from $w$.

\begin{lemma}\label{L:prime}
Any prime ideal $\maxim$ of $A$ containing $\pr$ is finitely generated and its extension $\maxim\ul A$ is again prime. 
\end{lemma}
\begin{proof}
Since $A/\pr A=\pol \kappa X$ is Noetherian, where $\kappa$ is the residue field of $\br $, the ideal $\maxim (A/\pr A)$ is finitely generated. Therefore so is $\maxim$, since by assumption $\pr$ is finitely generated. Moreover, $\ul A/\pr\ul A$ is the ultraproduct of the $\pol{\seq\kappa w}X$, so that by Theorem~\ref{T:SvdD}, the extension $\maxim (\ul A/\pr\ul A)$ is prime, whence so is $\maxim \ul A$.
\end{proof}

In particular, if $\seq\maxim w$ is an \sr\ of $\maxim$, then almost all $\seq\maxim w$ are prime ideals. Therefore, the following notions are well-defined (with the convention that we put $B_{\mathfrak n}$ equal to zero whenever $\mathfrak n$ is not a prime ideal of the ring $B$). Let $R$ be  a  local $\br$-affine algebra,  say, of the form $C_\maxim$, with $C$ a  finitely presented $\br$-algebra  and $\maxim$ a prime ideal  containing $\pr$.

\begin{definition}
 We call $(\ul C)_{\maxim \ul C}$ the \emph{non-standard hull} of $R$ and denote it $\ul R$. Moreover, if $\seq Cw$ and $\seq\maxim w$ are \sr{s} of $C$ and $\maxim$ respectively, then the collection $\seq Rw:=(\seq Cw)_{\seq\maxim w}$   is an \emph\sr\ of $R$.
\end{definition}

One easily checks that the ultraproduct of the \sr{s} $\seq Rw$ is precisely the non-standard hull $\ul R$.

\section{Flatness of Non-standard Hulls}

In this section, we specialize the notions from the previous result to the situation where each $\seq \br w$ is a \DVR. We fix throughout the following notation. For each $w$, let $\seq \br w$ be a \DVR\ with uniformizing parameter $\seq\pi w$ and with residue field $\seq\kappa w$. Let $\br$, $\pi$ and $\kappa$ be their respective ultraproducts, so that $\pi \br$ is the maximal ideal of $\br$ and $\kappa$ its residue field. The intersection of all $\pi^m\br$ is called the \emph{ideal of infinitesimals} of $\br$ and is denoted $\varpi$. Using \cite{SchEC}, one sees that $\br/\pi^m\br$ is an Artinian local Gorenstein $\kappa$-algebra of length $m$.

Fix   a finite tuple of variables $X$ and let $\ul A$ be the ultraproduct of the $\seq Aw:=\pol{\seq \br w}X$. Set $A:=\pol{\br}X$ and view it as a subring of $\ul A$.

\begin{proposition}\label{P:max}
For $I$ an ideal in $A$, the residue ring $A/I$ is Noetherian \iff\ $\varpi\subset I$. In particular, every maximal ideal of $A$ contains $\varpi$ and is of the form $\varpi A+J$ with $J$ a finitely generated ideal.
\end{proposition}
\begin{proof}
Let $C:=A/I$ for some ideal $I$ of $A$. If $C$ is Noetherian, then the intersection of all $\pi^nC$ is zero by Krull's Intersection Theorem. Hence $\varpi\subset I$. Conversely, if $\varpi\subset I$, then  since $A/\varpi A=\pol{(\br/\varpi)}X$ is Noetherian, so is $C$. The last assertion is now clear.
\end{proof}

In spite of Lemma~\ref{L:prime}, there are even maximal ideals of $A$ (necessarily not containing $\pi$) which do not extend to a proper ideal in $\ul A$. For instance with $X$ a single variable, the ideal $\varpi A+(1-\pi X)A$ is maximal (with residue field the field of fractions of $\br/\varpi$), but $\varpi\ul A+(1-\pi X)\ul A $ is the unit ideal. Indeed, let $\ul f$ be the ultraproduct of the 
	\begin{equation*}
	\seq fw:=(1-(\seq\pi w X)^w)/(1-\seq\pi w X).
	\end{equation*} 
Since $(1-\seq\pi w X)\seq fw\equiv 1$ modulo $(\seq\pi w)^w\seq Aw$, we get by \los\ that $(1-\pi X)\ul f\equiv 1$ modulo $\varpi\ul A$. Therefore, we cannot hope for $A\to \ul A$ to be faithfully flat. Nonetheless, using for instance a   result of \name{Aschenbrenner} on bounds of syzygies, we do have this property for local affine algebras. This result will prove to be crucial in what follows.

\begin{theorem}\label{T:ff}
The canonical \homo\ $A\to \ul A$ is flat. In particular, the canonical \homo\ of a  local $\br$-affine algebra to its non-standard hull is faithfully flat, whence in particular injective.
\end{theorem}
\begin{proof}
The last assertion is clear from the first, since the \homo\ $R\to \ul R$ is obtained as a   base change of $A\to \ul A$ followed by a suitable localization, for any  local $\br$-affine algebra $R$. I will provide two different proofs for the first assertion

For the first proof, we use a result of \name{Aschenbrenner} \cite{Asch} in order  to verify the equational criterion for flatness, that is to say, given a linear equation $L=0$, with $L$ a linear form over $A$, and given a solution $\ul{\mathbf f}$ over $\ul A$, we need to show that there exist solutions $\mathbf b_i$ in $A$ such that $\ul{\mathbf f}$ is an $\ul A$-linear combination of the $\mathbf b_i$. Choose $\seq Lw$ and $\seq{\mathbf f}w$   with respective ultraproducts $L$ and $\ul{\mathbf f}$.  In particular, almost all $\seq Lw$ have $\seq\br w$-complexity at most $c$, for some $c$ independent from $w$. By \los, $\seq{\mathbf f}w$ is a solution of the linear equation $\seq{L}w=0$, for almost all $w$. Therefore, by \cite[Corollary 4.27]{Asch}, there is a bound $c'$, only depending on $c$, such that  $\seq{\mathbf f}w$ is an  $\seq Aw$-linear combination of solutions $\seq{\mathbf b_1}w,\dots, \seq{\mathbf b_s}w$ of $\seq\br w$-complexity at most $c$. Note that $s$ can be chosen  independent from $w$ as well by \cite[Lemma 1]{SchBArt}. In particular, the ultraproduct $\mathbf b_i$ of the $\seq{\mathbf b_i}w$  lies in $A$ by Lemma~\ref{L:deg}. By \los, each $\mathbf b_i$ is a solution of $L=0$ in $\ul A$, whence in $A$, and $\ul{\mathbf f}$ is an $\ul A$-linear combination of the $\mathbf b_i$, proving flatness.

If we want to avoid the use of \name{Aschenbrenner}'s result, we can reason as follows. By Theorem~\ref{T:SvdD}, both extensions $A/\pi A\to \ul A/\pi \ul A$ and $A\tensor Q\to \ul A\tensor Q$ are faithfully flat, where $Q$ is the field of fractions of $\br$. Let $M$ be an $A$-module. Since $\pi$ is $A$-regular, the standard spectral sequence
 	\begin{equation*}
	 \tor {A/\pi A}p{\ul A/\pi\ul A}{\tor AqM {A/\pi A}} \implies \tor A{p+q}{\ul A/\pi\ul A}M
	\end{equation*}
 degenerates into short exact sequences
	\begin{multline*}
	\tor {A/\pi A}{i-1}{\ul A/\pi \ul A}{(0:_M\pi)} \to \tor Ai{\ul A/\pi\ul A}M \to \\
\tor {A/\pi A}i{\ul A/\pi\ul A}{M/\pi M},
	\end{multline*}
for all $i\geq 2$. In particular, the flatness of $A/\pi A\to \ul A/\pi\ul A$ implies the vanishing of  $\tor A2{\ul A/\pi\ul A}M$. Applying this to the short exact sequence
	\begin{equation*}
	\exactseq {\ul A} \pi {\ul A} {} {\ul A/\pi\ul A}
	\end{equation*}
we get a short exact sequence
	\begin{equation}\label{eq:unittor}
	0=\tor A2{\ul A/\pi\ul A}M\to \tor A1{\ul A}M \map \pi \tor A1{\ul A}M.
	\end{equation}
On the other hand, flatness of $A\tensor Q\to \ul A\tensor Q$ yields
	\begin{equation}\label{eq:torsiontor}
	\tor A1{\ul A}M\tensor Q = \tor{A\tensor Q}1{\ul A\tensor Q}{M\tensor Q}=0.
	\end{equation}

In order to prove flatness, it suffices by \cite[Theorem 7.8]{Mats} to show that $\tor A1{\ul A}{A/I}$ vanishes, for every finitely generated ideal $I$ of $A$. Towards a contradiction, suppose that $\tor A1{\ul A}{A/I}$ contains a non-zero element $\tau$. By \eqref{eq:torsiontor}, we have $a\tau=0$, for some non-zero $a\in\br$. As observed in \cite[Proposition 3]{Sab}, every polynomial ring over a valuation ring is coherent, so that in particular $I$ is finitely presented (namely,  since $I$ is torsion-free over $\br$, it is $\br$-flat, and therefore finitely presented by \cite[Theorem 3.4.6]{RayGru}). Hence we have some exact sequence
	\begin{equation*}
A^{a_2} \map{\varphi_2} A^{a_1} \map{\varphi_1}A\to A/I\to 0.
	\end{equation*}
Therefore $\tor A1{\ul A}{A/I}$ is calculated as the homology of the complex
	\begin{equation*}
	(\ul A)^{a_2} \map{\varphi_2} (\ul A)^{a_1} \map{\varphi_1} \ul A.
	\end{equation*}
Suppose $\tau$ is the image of a tuple $\mathbf x\in(\ul A)^{a_1}$ with $\varphi_1(\mathbf x)=0$. Hence $\mathbf x$ does not belong to $\varphi_2((\ul A)^{a_2})$ but $a\mathbf x$ does. Choose $\seq{\mathbf x}w$, $\seq aw$ and $\seq{\varphi_i}w$ with respective ultraproduct $\mathbf x$, $a$ and $\varphi_i$. By \los, almost all $\seq{\mathbf x}w$ lie in the kernel of $\seq{\varphi_1}w$ but not in the image of $\seq{\varphi_2}w$, yet $\seq aw\seq{\mathbf x}w$ lies in the image of $\seq{\varphi_2}w$. Choose $\seq nw\in\nat$ maximal such that $\seq{\mathbf y}w:=(\seq\pi w)^{\seq nw}\seq {\mathbf x}w$ does not lie in the image  of $\seq{\varphi_2}w$. Since almost all $\seq aw$ are non-zero, this maximum exists for almost all $w$. Therefore, if $\mathbf y$ is the ultraproduct of the $\seq{\mathbf y}w$, then $\varphi_1(\mathbf y)=0$ and $\mathbf y$ does not lie in $\varphi_2((\ul A)^{a_2})$, but $\pi\mathbf y$ lies in $\varphi_2((\ul A)^{a_2})$. Therefore, the image of $\mathbf y$ in $\tor A1{\ul A}{A/I}$ is a non-zero element annihilated by $\pi$, contradicting \eqref{eq:unittor}.
\end{proof}

\begin{remark}
In \cite{SchUBS}, I exhibit a general connection  between the flatness of an ultraproduct over certain canonical subrings  and the existence of bounds on syzygies. In particular, using these ideas, the second argument in the above proof of flatness reproves the result in \cite{Asch}. In fact,   the role played here by coherence  is not accidental either; see \cite{AschSab} or \cite{SchUBS} for more details.
\end{remark}

\begin{theorem}\label{T:dom}
Let $R$ be a local $\br$-affine algebra with non-standard hull $\ul R$ and \sr\ $\seq Rw$. 
\begin{itemize}
\item Almost all $\seq Rw$ are flat over $\seq\br w$ \iff\ $R$ is torsion-free over $\br$ \iff\ $\pi$ is $R$-regular.
\item Almost all $\seq Rw$ are domains \iff\ $R$ is.
\end{itemize}
 \end{theorem}
\begin{proof}
Suppose first that almost all $\seq Rw$ are flat over $\seq\br w$, which amounts in this case, to almost all $\seq Rw$ being torsion-free over $\seq\br w$. By \los, $\ul R$ is torsion-free over $\br$, and since $R\subset \ul R$, so is $R$. Conversely, assume $\pi$ is $R$-regular. By faithful flatness, $\pi$ is $\ul R$-regular, whence almost all  $\seq\pi w$ are $\seq Rw$-regular by \los. Since the $\seq\br w$ are \DVR{s}, this means that almost all $\seq\br w\to \seq Rw$ are flat.

If almost all $\seq Rw$ are domains, then so is $\ul R$ by \los, and hence so is $R$, since it embeds in $\ul R$. Conversely, assume $R$ is a domain. If $\pi=0$ in $R$, then $\ul R$ is a domain by Lemma~\ref{L:prime}, whence so are almost all $\seq Rw$ by \los. So assume $\pi$ is non-zero in $R$, whence $R$-regular. By what we just proved, $R$ is then torsion-free over $\br$. Let $Q$ be the field of fractions of  $\br$. Write $R$ in the form $S/\pr$, where $S$ is some localization of $A$ at a prime ideal containing $\pi$ and $\pr$ is a finitely generated prime ideal in $S$. Since $S/\pr$ is torsion-free over $\br$, the extension $\pr(S\tensor_\br Q)$ is again prime and its contraction in $S$ is $\pr$. By Theorem~\ref{T:SvdD}, since we are now over a field, $\pr(\ul S\tensor_\br Q)$ is a prime ideal, where $\ul S$ is the non-standard hull of $S$ (note that $\ul S\tensor_\br Q$ is then the non-standard hull of $S\tensor_\br Q$ in the sense of \cite{SchNSTC}). Moreover, since $S/\pr$ is torsion-free over $\br$, so is $\ul S/\pr\ul S$ by the first assertion. This in turn means that 
	\begin{equation*}
	\pr\ul S= \pr (\ul S\tensor_\br Q)\cap \ul S,
	\end{equation*} 
showing that $\pr\ul S$ is prime. It follows then from \los\ that  almost all $\seq\pr w$ are prime, where $\seq\pr w$ is an \sr\ of $\pr$, and hence almost all $\seq Rw$ are domains.
\end{proof}

\begin{remark}\label{R:extprim}
The last assertion is equivalent with saying that any prime ideal in $R$ extends to a prime ideal in $\ul R$. Indeed, let $\mathfrak q$ be a prime ideal in $R$ with \sr\  $\seq{\mathfrak q} w$. By the above result (applied to $R/{\mathfrak q}$ and its \sr\ $\seq Rw/\seq {\mathfrak q}w$), we get that almost all $\seq{\mathfrak q} w$ are prime, whence so is their ultraproduct ${\mathfrak q}\ul R$, by \los.  
\end{remark}

\section{\Qdim}\label{s:qdim}

In this and the next section, we will study the local algebra of the category  $\affine\br$. Although part of the theory can be developed for arbitrary base rings $\br$, we will only deal with the case that $\br$ is a local domain of embedding dimension one. Recall that the \emph{embedding dimension} of a local ring $(Z,\pr)$ is by definition the minimal number of generators of $\pr$, and its \emph{ideal of infinitesimal} $\varpi$ is the intersection of all powers $\pr^n$. Of course, if $Z$ is moreover Noetherian, then its ideal of infinitesimals is zero.  In general, we call $Z/\varpi$ the \emph\sepquot\ of $Z$.

For the duration of the next two sections, let $\br$  denote a local domain of embedding dimension one, with generator of the maximal ideal $\pi$, with ideal of infinitesimals $\varpi$ and with residue field $\kappa$.

\begin{lemma}\label{L:dvr}
The \sepquot\  $\br/\varpi$ of $\br$ is a \DVR\ with uniformizing parameter $\pi$.
\end{lemma}
\begin{proof}
For each element $a\in\br$ outside $\varpi$, there is a smallest $e\in\nat$ for which $a\notin\pi^{e+1}\br$. Hence  $a=u\pi^e$ with $u$ a unit in $R$. It is now straightforward to check that the assignment $a\mapsto e$ induces a discrete valuation on $\br/\varpi$. 
\end{proof}

Note that we do not even need $\br$ to be domain, having positive depth (that is to say, assuming that $\pi\br$ is not an associated prime of $\br$; see \cite[Proposition 9.1.4]{BH}) would suffice, for then $\pi$ is necessarily $\br$-regular. However, we do not need this amount of generality as in all our applications $\br$ will be of the following special form.

\begin{definition}\label{D:uDVR}
We say that $\br$ is an \emph{ultra-DVR}, if it is realized as the ultraproduct of some \DVR{s} $\seq \br w$.
\end{definition} 

 Note that $\pi$ and $\kappa$ are then the respective ultraproducts of the uniformizing parameter  $\seq\pi w$ and the residue field  $\seq\kappa w$ of $\seq \br w$.

We will work in the category $\affine\br$ of local $\br$-affine algebras (see the Introduction), that is to say, the category of algebras of the form $(A/I)_\maxim$, with $I$ a finitely generated ideal and $\maxim$ a prime ideal containing $\pi$ (and $I$). Nonetheless, some results can be stated even for local algebras which are locally finitely generated over $\br$, that is without the assumption that $I$ is finitely generated. We call $R$ a \emph{torsion-free} $\br$-algebra if it is torsion-free over $\br$ (that is to say, if $ar=0$ for some $r\in R$ and some non-zero $a\in \br$, then $r=0$). Recall from Theorem~\ref{T:dom} that  a local $\br$-affine algebra $R$ is torsion-free \iff\ $\pi$ is $R$-regular.

If   $\br$ is  an ultra-DVR, then we set $\seq Aw:=\pol{\seq \br w}X$ and let $\ul A$ be their ultraproduct. If $R$ belongs to   $\affine\br$, that is to say, is a local $\br$-affine algebra of the form $(A/I)_\maxim$, then $\ul R:=(\ul A/I\ul A)_{\maxim\ul A}$ denotes its  non-standard hull and we let $\seq Rw:=(\seq Aw/\seq Iw)_{\seq\maxim w}$ be an \sr\ of $R$, where $\seq Iw$ and $\seq \maxim w$ are \sr{s} of $I$ and $\maxim$ respectively. Note that $\maxim$ is finitely generated, as it contains  by definition $\pi$.

\begin{lemma}\label{L:inf}
Let $R$ be a local ring which is locally finitely generated over $\br$. If $I$ is a proper ideal in $R$ containing some power $\pi^m$, then the intersection of all $I^N$ is equal to $\varpi R$. 
\end{lemma}
\begin{proof}
 Suppose $\pi^m\in I\subset\maxim$. Let $J$ be the  intersection of all $I^N$. Since $\pi^m\in I$, we get that  $\varpi R\subset J$. Since  $S:=R/\varpi R$ is locally finitely generated over the \DVR\ $\br/\varpi$ (see Lemma~\ref{L:dvr}), it is itself Noetherian. Applying Krull's Intersection Theorem (see for instance \cite[Theorem 8.10]{Mats}), we get that $JS=(0)$, and hence that $J= \varpi R$. 
\end{proof}

\begin{lemma}\label{L:comp}
Let $\br$ be  an ultra-DVR (see Definition~\ref{D:uDVR}) and let $(R,\maxim)$ be a local $\br$-affine algebra with \sr\ $\seq Rw$.  The $\maxim$-adic completion $\complet R$ of $R$ is isomorphic to $\ul R/\mathfrak i$, where $\mathfrak i$ is the ideal of infinitesimals of $\ul R$. 

In particular, $\complet R$ is Noetherian.
\end{lemma}
\begin{proof}
Clearly, $R$ is $\maxim$-adically dense in $\ul R$, so that they have a common completion. On the other hand, by saturatedness of ultraproducts, $\ul R$ is quasi-complete in the sense that every Cauchy sequence has a (non-unique) limit. Therefore, its \sepquot\ $\ul R/\mathfrak i$ is complete.
\end{proof}

Our first goal is to introduce a good notion of dimension. Below, the \emph{dimension} of a ring will always mean its \emph{Krull dimension}. Recall that it is always finite for Noetherian local rings.

\begin{theorem}\label{T:qdim}
For  a local ring $(R,\maxim)$ which is locally finitely generated over $\br$, the following numbers are all equal: 
\begin{itemize}
\item the minimal length $d$ of a tuple in $R$ generating an $\maxim$-primary ideal;
\item the dimension $\complet d$ of the completion $\complet R$;
\item the dimension $\tilde d$ of $R/\varpi R$;
\item the degree $\underline d$ of the \emph{Hilbert-Samuel polynomial} $\chi_R$, where $\chi_R$ is the unique polynomial  with rational coefficients for which $\chi_R(n)$ equals the length of $R/\maxim^{n+1}$ for all large $n$.
\end{itemize}
Moreover, if $\pi$ is $R$-regular, then  $R/\pi R$ has dimension $d-1$. 

If $\br$ is  an ultra-DVR and $R$ is a torsion-free local $\br$-affine algebra  with \sr\ $\seq Rw$, then almost all $\seq Rw$ have dimension   $d$.
\end{theorem}
\begin{proof}
Let $\tilde R:=R/\varpi R$, so that $\tilde R$ is a Noetherian local ring (see  the proof of Lemma~\ref{L:inf}) and $\tilde d$ is finite. Since $R/\maxim^n\iso \tilde R/\maxim^n\tilde R$, the completion of $\tilde R$ is  $\complet R$ and hence, $\complet d=\tilde d$. Moreover, $\chi_R=\chi_{\complet R}$, so that by the Hilbert-Samuel theory, $\underline d=\complet d$.

Let $\mathbf x$ be a tuple of length $\tilde d$ such that its image in $\tilde R$ is a system of parameters of $\tilde R$. Hence, for some $n$, we have that $\maxim^n\subset \mathbf xR+\varpi R$. In particular, since $\varpi R\subset \pi^{n+1}R$ by Lemma~\ref{L:inf}, we can find $x\in\mathbf xR$ and $r\in R$, such that $\pi^n=x+r\pi^{n+1}$. Therefore,   $\pi^n\in\mathbf xR$, since $(1-r\pi)$ is a unit. Since $\varpi\subset\pi^n\br$, we get that $\maxim^n\subset \mathbf xR$, showing that $\mathbf xR$ is an $\maxim$-primary ideal and hence that $d\leq \tilde d$. On the other hand, if  $\mathbf y$ is a tuple of length $d$ such that $\mathbf yR$ is $\maxim$-primary, then $\mathbf y\tilde R$ is an $\maxim \tilde R$-primary ideal, and hence $\tilde d\leq d$. This concludes the proof of the first assertion.

Assume that $\pi$ is moreover $R$-regular. I claim that $\pi$ is $\tilde R$-regular. Indeed, suppose $\pi \tilde r=0$, for some $\tilde r\in \tilde R$. Take a pre-image $r\in R$, so that $\pi r\in\varpi R\subset\pi^nR$, for every $n$. Since $\pi$ is $R$-regular, we get that $r\in\pi^{n-1}R$, for all $n$. Therefore $r\in\varpi R$, whence $\tilde r=0$ in $\tilde R$, as we needed to show. Since $\pi$ is $\tilde R$-regular and $\tilde R/\pi\tilde R=R/\pi R$, the dimension of  $R/\pi R$ is $\tilde d-1$.

Suppose finally that $\br$ is moreover an ultra-DVR. We already observed that $\seq Rw/\seq\pi w\seq Rw$ is an \sr\ of $R/\pi R$ in the sense of \cite{SchNSTC}. In particular, by  \cite[Theorem 4.5]{SchNSTC}, almost all $\seq Rw/\seq\pi w\seq Rw$ have dimension $\tilde d-1$. Since $\pi$ is $\ul R$-regular by flatness, whence $\seq\pi w$ is $\seq Rw$-regular by \los, we get that $\seq Rw$ has dimension $\tilde d$, for almost all $w$.
\end{proof}

\begin{definition}[\Qdim]\label{D:qdim}
The minimal number of generators of an $\maxim$-primary ideal is called the  \emph\qdim\ of $R$. We call a tuple $\mathbf x$ in $R$ \emph{generic}, if it generates an $\maxim$-primary ideal and has length equal to the \qdim\ of $R$.
\end{definition}

\begin{corollary}\label{C:trim}
In   a local ring $(R,\maxim)$ which is locally finitely generated over $\br$, every tuple generating an $\maxim$-primary ideal can be trimmed to a generic sequence by omitting some of its entries.
 \end{corollary}
\begin{proof}
Let $\tilde R:=R/\varpi R$ and let $d$ be the \qdim\ of $R$. Let $\mathbf x$ be a tuple generating an $\maxim$-primary ideal of $R$. By a well-known prime avoidance argument in the Noetherian ring $\tilde R$, we can rearrange this tuple so that it has the form $(\mathbf y,\mathbf z)$ with $\mathbf y$ a system of parameters in $\tilde R$. In particular, $\norm{\mathbf y}=d$ by Theorem~\ref{T:qdim}. Let $S:=R/\mathbf yR$ and $\tilde S:=S/\varpi S$. By Theorem~\ref{T:qdim}, the \qdim\ of $S$ is equal to the dimension of $\tilde S$, whence is zero since $\tilde S=\tilde R/\mathbf y\tilde R$. Therefore, the empty tuple is a generic sequence in $S$, and hence   $\mathbf yR$ is $\maxim$-primary. Since $\mathbf y$ has length equal to the \qdim\ of $R$, it is therefore a generic sequence.
\end{proof}

In fact the above proof shows that there is a one-one correspondence between generic sequences in $R$ and systems of parameters in $R/\varpi R$. In general, the last assertion in Theorem~\ref{T:qdim} is false when $R$ is not torsion-free. For instance,  let $R:=\br/\gamma\br$ with $\gamma$ a non-zero infinitesimal, so that each $\seq Rw=\seq\br w/\seq\gamma w\seq\br w$ has dimension zero, but $R/\varpi R$ is the (one-dimensional) \DVR\ $\br/\varpi$. 

In the following definition, let $\br$ be  an ultra-DVR and let $R$ be a local $\br$-affine algebra of \qdim\ $d$, with \sr\ $\seq Rw$. Note that the $\seq Rw$ have almost all  dimension at most $d$. Indeed, if $\mathbf y$ has length $d$ and generates an $\maxim$-primary ideal,  then almost all $\seq{\mathbf y}w$ are $\seq\maxim w$-primary by \los, for $\seq{\mathbf y}w$ an \sr\ of $\mathbf y$.

\begin{definition}
We say that $R$ is \emph\serene\ if almost all $\seq Rw$ have dimension equal to the \qdim\ of $R$.
\end{definition}

Theorem~\ref{T:qdim} shows that every  torsion-free local $\br$-affine algebra is \serene. In particular, over an ultra-DVR, the \blim\ $R$ of domains $\seq Rw$ of uniformly bounded $\seq\br w$-complexity is \serene, since $\ul R$ is then a domain by \los, whence so is $R$ as it embeds in $\ul R$. The next result shows that generic sequences in an \serene\ ring are the analog of systems of parameters.

\begin{corollary}\label{C:gensop}
Let $\br$ be an ultra-DVR and $R$ an \serene\  local $\br$-affine algebra with \sr\ $\seq Rw$.   Let $\mathbf x$ be a tuple in $R$  with \sr\ $\seq{\mathbf x}w$. 

If $\mathbf x$ is generic, then $\seq{\mathbf x}w$ is a system of parameters of $\seq Rw$, for almost all $w$. Conversely,  if $(\seq\pi w)^c\in \seq{\mathbf x}w\seq Rw$, for some $c$ and almost all $w$, then $\mathbf x$ is generic.
\end{corollary}
\begin{proof}
Let $\maxim$ be the maximal ideal of $R$, with \sr\ $\seq\maxim w$. Let  $d$ be the \qdim\ of $R$, so that almost all $\seq Rw$ have dimension $d$. Suppose first that $\mathbf x$ is generic, so that $\norm{\mathbf x}=d$ and  $\mathbf xR$ is $\maxim$-primary. Since $\mathbf x\ul R$ is then $\maxim\ul R$-primary, $\seq{\mathbf x}w\seq Rw$ is $\seq\maxim w$-primary  by \los, showing that $\seq{\mathbf x}w$ is a system of parameters for almost all $w$. 

Conversely, suppose  $\seq{\mathbf x}w$ is a system of parameters of $\seq Rw$, generating an ideal containing  $(\seq\pi w)^c$. If $\mathbf x$ is the ultraproduct of the $\seq{\mathbf x}w$, then  it is defined over $R$ by Lemma~\ref{L:deg}. By \los\ and faithful flatness, $\pi^c\in\mathbf xR$. Applying \cite[Corollary 4]{SchBArt} to the Artinian base ring $\seq \br w/(\seq\pi w)^c$, we can find  a bound $c'$, only depending on $c$, such that $(\seq\maxim w)^{c'} \subset  \seq{\mathbf x}w\seq Rw$, for almost all $w$. Hence $\maxim^{c'}\ul R\subset\mathbf x\ul R$, so that by faithful flatness, $\mathbf xR$ is $\maxim$-primary. This shows that $\mathbf x$ is generic.
\end{proof}

The additional requirement in the converse is necessary: indeed, for arbitrary $\seq nw>0$, the element $(\seq\pi w)^{\seq nw}$  is a parameter in $\seq \br w$ and has $\seq \br w$-complexity zero, but if $\seq nw$ is unbounded, its ultraproduct is an infinitesimal whence not generic. To characterize \serene\ rings, we use the following notion introduced in \cite{SchABCM}. 

\begin{definition}[Parameter Degree]
The \emph{parameter degree} of a Noetherian local ring $C$ is by definition the smallest possible length of a residue ring $C/\mathbf xC$, where $\mathbf x$ runs over all systems of parameters of $C$.
\end{definition}

In general, the parameter degree is larger than the multiplicity, with equality precisely when $C$ is \CM\ (see \cite[Theorem 17.11]{Mats}). The homological degree of $C$ is an upper bound for its parameter degree (see \cite[Corollary 4.6]{SchABCM}). A priori, being \serene\ is a property of the \sr{s} of $R$, of for that matter, of its non-standard hull. However, the last equivalent condition in the next result shows that it is in fact an intrinsic property

\begin{proposition}\label{P:pardeg}
Let $\br$ be  an ultra-DVR and let $R$ be a local $\br$-affine algebra with \sr\ $\seq Rw$. The following are equivalent:
\begin{enumerate}
\item\label{i:eqdim} $R$ is \serene;
\item\label{i:ram} there exists a $c\in\nat$, such that for almost all $w$, we can find a system of parameters $\seq{\mathbf x}w$ of $\seq Rw$ of $\seq\br w$-complexity at most $c$, generating an ideal containing  $(\seq\pi w)^c$;
\item\label{i:pardeg} there exists an $e\in\nat$, such that almost all $\seq Rw$ have parameter degree at most $e$;
\item\label{i:intr} for every generic sequence in $R$ of the form $(\pi,\mathbf y)$, we have that $\mathbf yR\cap \br$ is the zero ideal.
\end{enumerate}
\end{proposition}
\begin{proof}
Let $\maxim$ be the maximal ideal of $R$, with \sr\ $\seq\maxim w$. Let $d$ be the \qdim\ of $R$ and let $d'$ be the dimension of almost all $\seq Rw$. Suppose first that $d=d'$. Let $\mathbf x$ be any generic sequence in $R$ with \sr\ $\seq{\mathbf x}w$. By \los, almost all $\seq{\mathbf x}w$ generate an $\seq\maxim w$-primary ideal. Since their length is equal to the dimension of $\seq Rw$, they are almost all systems of parameters of $\seq Rw$. Choose $c$ large enough so that $\pi^c\in\mathbf xR$. Enlarging $c$ if necessary, we may moreover assume by Lemma~\ref{L:deg} that almost all $\seq{\mathbf x}w$ have $\seq\br w$-complexity at most $c$. By \los, $(\seq\pi w)^c\in \seq{\mathbf x}w\seq Rw$, so that \eqref{i:ram} holds.

Assume next that $c$ and the $\seq{\mathbf x}w$ are as in \eqref{i:ram}. Let $\seq{\overline R}w:=\seq Rw/(\seq \pi w)^c\seq Rw$. We can apply \cite[Corollary 2]{SchBArt} over $\seq\br w/(\seq\pi w)^c\seq\br w$ to the $\seq\maxim w\seq{\overline R}w$-primary ideal $\seq{\mathbf x}w\seq {\overline R}w$, to conclude that there is some $c'$, depending only on $c$, such that $\seq {\overline R}w/\seq{\mathbf x}w\seq {\overline R}w$ has length at most $c'$. Since the latter residue ring is just $\seq Rw/\seq{\mathbf x}w\seq Rw$ by assumption, the parameter degree of $\seq Rw$ is at most $c'$, and hence \eqref{i:pardeg} holds.

To show that \eqref{i:pardeg} implies \eqref{i:eqdim}, assume that almost all $\seq Rw$ have parameter degree at most $e$. Let $\seq{\mathbf y}w$ be a system of parameters of $\seq Rw$ such that $\seq Rw/\seq{\mathbf y}w\seq Rw$ has length at most $e$, for almost all $w$. It follows that $(\seq\maxim w)^e$ is contained in $\seq{\mathbf y}w\seq Rw$. Let $\ul{\mathbf y}$ be the ultraproduct of the $\seq{\mathbf y}w$. By \los, $\maxim^e\ul R\subset \ul{\mathbf y}\ul R$. Since the completion $\complet R$ is a homomorphic image of $\ul R$ by Lemma~\ref{L:comp}, we get that $\ul{\mathbf y}\complet R$ is $\maxim\complet R$-primary. Since $\ul{\mathbf y}$ has length at most $d'$ (some entries might be zero in $\complet R$), the dimension of $\complet R$ is at most $d'$. Since we already remarked that $d'\leq d$, we get from Theorem~\ref{T:dom} that $d'=d$.

So remains to show that \eqref{i:intr} is equivalent to the other conditions. Assume first that it holds but that $R$ is not \serene. Since we have inequalities $d-1\leq d'\leq d$, this means that $d'=d-1$. Moreover, $R/\pi R$ must have \qdim\ also equal to $d-1$, for if not, its \qdim\ would be $d$, whence almost all $\seq Rw/\seq\pi w\seq Rw$ would have dimension $d$ by \cite[Theorem 4.5]{SchNSTC}, which is impossible. Since there is a uniform bound $c$ on the $\seq\br w$-complexity of each $\seq Rw$, we can choose using Remark~\ref{R:pa}, a system of parameters $\seq{\mathbf y}w$ of $\seq Rw$ of $\seq\br w$-complexity at most $c$ (and hence of length $d-1$). In particular, some power of $\seq\pi w$ lies in $\seq{\mathbf y}w\seq Rw$. Let $a\in\br$ be the ultraproduct of these powers. If $\mathbf y$ is the ultraproduct of the $\seq{\mathbf y}w$, then $\mathbf y$ is already defined over $R$ by Lemma~\ref{L:deg}. By \los, $a\in\mathbf y\ul R$, whence by faithful flatness, $a$ is a non-zero element in $\mathbf yR\cap\br$. Therefore, to reach the desired contradiction with \eqref{i:intr}, we only need to show that $(\pi,\mathbf y)$ is generic. As we already established, $\seq Rw/\seq\pi w\seq Rw$ has dimension $d-1$, so that $\seq{\mathbf y}w$ is also a system of parameters in that ring. Therefore, $\mathbf y$ is a system of parameters in $R/\pi R$ by \cite[Theorem 4.5]{SchNSTC}. This in turn implies that $(\pi,\mathbf y)$ generates an $\maxim$-primary ideal in $R$. Since this tuple has length $d$, it is therefore generic, as we wanted to show.

Finally, assume $R$ is \serene, and suppose $(\pi,\mathbf y)$ is generic. Let $a\in\mathbf yR\cap \br$ and choose \sr{s} $\seq aw$ and $\seq{\mathbf y}w$ of $a$ and $\mathbf y$ respectively. By \los, $\seq aw\in\seq{\mathbf y}w\seq Rw$. However, if $a$ is non-zero, then $\seq aw$ is a power of $\seq\pi w$, which contradicts the fact that $(\seq\pi w,\seq{\mathbf y}w)$ is a system of parameters by Corollary~\ref{C:gensop}. So $a=0$, as we needed to show.
\end{proof}

\begin{corollary}\label{C:bdpd}
For each $c$, there exists a bound $\bound{PD}c$ with the following property. Let $V$ be a \DVR\ and let $R$ be a local $V$-algebra of $V$-complexity at most $c$. If $R$ is torsion-free over $V$, then the parameter degree of $R$ is at most $\bound{PD}c$.
\end{corollary}
\begin{proof}
If the statement is false for some $c$, then we can find for each $w$ a \DVR\ $\seq\br w$ and a torsion-free local $\seq \br w$-algebra $\seq Rw$ of $\seq\br w$-complexity at most $c$, whose parameter degree is at least $w$. Let $R$ be the \blim\ of the $\seq Rw$ and let $\ul R$ be their ultraproduct. Since $\seq\pi w$ is $\seq Rw$-regular, $\pi$ is $\ul R$-regular, whence $R$-regular. Hence $R$ is \serene\ by Theorem~\ref{T:qdim}. Therefore, there is a bound on the parameter degree of almost all $\seq Rw$ by Proposition~\ref{P:pardeg}, contradicting our assumption.
\end{proof}

 Our next goal is to introduce a notion similar to height. Let $I$ be an arbitrary ideal of $R$.

\begin{definition}[\Qht]
We call the \emph\qht\ of $I$ the maximum of all $h$ such that there exists a generic sequence with its first $h$ entries in $I$.
\end{definition}

For Noetherian rings, we cannot expect a good relationship between the height of an ideal and the dimension of its residue ring, unless the ring is a catenary domain; the following is a \q\ analogue.

\begin{theorem}\label{T:qht}
Let $\br$ be  an ultra-DVR and let $R$ be a local $\br$-affine domain with \sr\ $\seq Rw$. Let $I$ be a finitely generated ideal in $R$ with \sr\ $\seq Iw$. 

If $R/I$ is \serene, then the \qht\ of $I$ is equal to the \qdim\ of $R$ minus the \qdim\ of $R/I$, and this is also equal to the height of almost all $\seq Iw$.
\end{theorem}
\begin{proof}
Let $d$ be the \qdim\ of $R$ and $\overline d$ the \qdim\ of $R/I$. Since  a domain is \serene, almost all $\seq Rw$ have dimension $d$ by Theorem~\ref{T:qdim}, and by assumption, almost all $\seq Rw/\seq Iw$ have \qdim\ $\overline d$. Let $h$ be the \qht\ of $I$. Let $\mathbf z$ be a generic sequence in $R$ with its first $h$ entries in $I$, and let $\seq{\mathbf z}w$ be an \sr\ of $\mathbf z$. By Corollary~\ref{C:gensop}, almost all $\seq{\mathbf z}w$ are a system of parameters in $\seq Rw$. Since by \los\ the first $h$ entries of $\seq{\mathbf z}w$ lie in $\seq Iw$, we get that $\seq Rw/\seq Iw$ has dimension at most $d-h$. In other words, $h\leq d-\overline d$. Since almost all $\seq Rw$ are catenary domains,  almost all $\seq Iw$ have height $d-\overline d$. 

Let $\mathbf x$ be a  tuple of length $\overline d$ whose image in $R/I$ is generic. Put $S:=R/\mathbf xR$ and let $e$ be its \qdim. If $\seq{\mathbf x}w$ is an \sr\ of $\mathbf x$, then almost all $\seq{\mathbf x}w$ are a system of parameters in $(\seq Rw/\seq Iw)$, by Corollary~\ref{C:gensop}. Since $\seq{\mathbf x}w$ is therefore part of a system of parameters of $\seq Rw$, we get that $\seq Sw:=\seq Rw/\seq{\mathbf x}w\seq Rw$ has dimension $d-\overline d$ by \cite[Theorem 14.1]{Mats}.  Using Corollary~\ref{C:trim}, we can find a tuple $\mathbf y$  of length $e$ with entries in $I$, so that its image in $S$ is a generic sequence. It follows that $\mathbf xR+\mathbf yR$ is $\maxim$-primary. We already observed that there are only two possibilities for $e$, to wit,  $d-\overline d$ or $d-\overline d+1$. If $e=d-\overline d$, then $(\mathbf x,\mathbf y)$ has length $d$ and hence is generic. Since $\mathbf y$ has all its entries in $I$, we get that $h$ is at least $d-\overline d$ and we are done in this case. 

Hence assume that $e=d-\overline d+1$, so that $(\mathbf x,\mathbf y)$ has length $d+1$.  By Corollary~\ref{C:trim}, we can trim  this tuple to a generic sequence, by omitting an appropriately chosen entry. Therefore, $I$ contains at least $e-1$ entries of this trimmed sequence, so that $h\geq e-1=d-\overline d$ and again we are done.
\end{proof}

\section{\Q\ Singularities}\label{s:qsing}

In this section, we maintain the notation introduced in the previous section. Our goal is to extend several singularity notions of Noetherian local rings to the category of local $\br$-affine algebras.

\subsubsection*{Grade and Depth.}
Let $B$ be an arbitrary ring and $I:=\rij xnB$ a finitely generated ideal. The \emph{grade} of $I$, denoted $\op{grade}(I)$, is by definition equal to $n-h$, where $h$ is the largest value $i$ for which the $i$-th Koszul homology $H_i\rij xn$ is non-zero. For   a local ring $R$ of finite embedding dimension, we define its \emph{depth} as the grade of its maximal ideal. 

If $B$ is moreover Noetherian, then we can define the grade of $I$ alternatively as the minimal $i$ for which $\ext Bi{B/I}B$ is non-zero (for all this see for instance \cite[\S9.1]{BH}). An arbitrary local ring has positive depth \iff\ its maximal ideal is not an associated prime.   Grade, and hence depth, \emph{deforms well}, in the sense that the 
	\begin{equation*}
	\op{grade}(I(B/\mathbf xB))= \op{grade}(I) -\norm{\mathbf x}
	\end{equation*}
for every $B$-regular sequence $\mathbf x$. For a  locally finitely generated $\br$-algebra $(R,\maxim)$, its depth never exceeds its \qdim. Indeed, by definition, the grade of a finitely generated ideal never exceeds its minimal number of generators, and by \cite[Proposition 9.1.3]{BH}, the depth of $R$ is equal to the grade of any of its $\maxim$-primary ideals. It follows that the depth of $R$ is at most its \qdim. 

In general,  the grade of a finitely generated ideal might be positive without it containing a $B$-regular element.   However, the next lemma shows that this is not the case for ultraproducts of Noetherian local rings.

\begin{lemma}\label{L:nd}
Let $\ul C$ be the ultraproduct of Noetherian local rings $\seq Cw$ and let $\ul I$ be a finitely generated ideal of $\ul C$ obtained as the ultraproduct of ideals $\seq Iw\subset\seq Cw$. 

If $\ul I$ has grade $n$, then there exists a $\ul C$-regular sequence of length $n$ with all of its entries in $\ul I$. Moreover, any permutation of a $\ul C$-regular sequence is again $\ul C$-regular.
\end{lemma}
\begin{proof}
By \cite[Proposition 9.1.3]{BH}, there exists a finite tuple of variables $Y$ and a $\pol{\ul C}Y$-regular sequence $\ul{\mathbf f}$ of length $n$, with all of its entries in $\ul I\pol {\ul C}Y$. Choose  tuples $\seq{\mathbf f}w$ in $\pol {\seq Cw}Y$ so that their ultraproduct is $\ul{\mathbf f}$. By \los, $\seq{\mathbf f}w$ is $\pol {\seq Cw}Y$-regular and has all of its entries in $\seq Iw\pol {\seq Cw}Y$, for almost all $w$. This shows that $\seq Iw\pol{\seq Cw}Y$ has grade at least $n$. Since $\seq Cw\to \pol{\seq Cw}Y$ is faithfully flat, $\seq Iw$ has grade at least $n$ by \cite[Proposition 9.1.2]{BH}. Hence, since $\seq Cw$ is Noetherian, we can find a $\seq Cw$-regular sequence $\seq{\mathbf x}w$ of length $n$   with all of its entries in $\seq Iw$. By \los, the ultraproduct $\ul{\mathbf x}$   of the $\seq{\mathbf x}w$  is $\ul C$-regular and has all of its entries in $\ul I$. 

The last assertion follows from \los\ and the fact that in a Noetherian local ring, any permutation of a regular sequence is again regular (\cite[Theorem 16.3]{Mats}).
\end{proof}

Recall that a Noetherian local ring for which its dimension and its depth (respectively, its dimension and its embedding dimension) coincide is \CM\ (respectively, regular). We will shortly see that upon replacing dimension by \qdim, we get equally well behaved notions. Let us therefore make the following definitions, for $R$ a locally finitely generated $\br$-algebra.

\begin{definition}
We say that $R$ is \emph\qCM, if its \qdim\ is equal to its depth, and \emph\qreg, if its \qdim\ is equal to its embedding dimension.
\end{definition}

Clearly, a \qreg\ local ring is \qCM.

\begin{theorem}\label{T:qCM}
Let $\br$ be  an ultra-DVR and let $R$ be an \serene\  local $\br$-affine algebra with \sr\ $\seq Rw$. In order for  $R$ to be \qCM\ it is necessary and sufficient that almost all $\seq Rw$ are \CM.
\end{theorem}
\begin{proof}
Let $d$ be the \qdim\ of $R$ and $\delta$ its depth. Suppose first that $d=\delta$. Since $R\to \ul R$ is faithfully flat, $\ul R$ has depth $\delta$ as well by \cite[Proposition 9.1.2]{BH}. By Lemma~\ref{L:nd}, there is an $\ul R$-regular sequence $\ul{\mathbf x}$ of length $d$. If $\seq{\mathbf x}w$ is an \sr\ of $\ul{\mathbf x}$, then almost all $\seq{\mathbf x}w$ are $\seq Rw$-regular by \los. Since almost all $\seq Rw$ have dimension $d$ by assumption, we showed that they are \CM.

Conversely, assume almost all $\seq Rw$ are \CM. It follows by reversing the above argument that $\ul R$ has depth $d$ and hence, so has $R$, by faithful flatness.
\end{proof}

Since every system of parameters is a regular sequence in a local \CM\ ring, we expect a similar behavior for generic sequences, and this indeed holds.

\begin{theorem}\label{T:genreg}
Let $\br$ be an ultra-DVR and let $R$ be an \serene\  local $\br$-affine algebra. If $R$ is \qCM, then any generic sequence is $R$-regular.
\end{theorem}
\begin{proof}
Let $\mathbf x$ be a generic sequence with \sr\ $\seq{\mathbf x}w$. Almost all $\seq{\mathbf x}w$ are a system of parameters in $\seq Rw$, by Corollary~\ref{C:gensop}. Since almost all $\seq Rw$ are \CM\ by Theorem~\ref{T:qCM}, almost all $\seq{\mathbf x}w$ are $\seq Rw$-regular. Hence $\mathbf x$ is $\ul R$-regular, by \los, whence $R$-regular, by faithful flatness.
\end{proof}

\begin{theorem}\label{T:qreg}
Let $\br$ be an ultra-DVR. An \serene\  local $\br$-affine algebra $R$ with \sr\ $\seq Rw$, is \qreg\ \iff\ almost all $\seq Rw$ are regular local rings.
\end{theorem}
\begin{proof}
Let $\maxim$ be the maximal ideal of $R$, with \sr\ $\seq\maxim w$. Let $\ul R$ be the non-standard hull of $R$. Let $\epsilon$ be the embedding dimension of $R$ and $d$ its \qdim. Suppose that $R$ is \qreg, that is to say, that $\epsilon=d$. Hence $\maxim=\mathbf xR$ for some $d$-tuple $\mathbf x$ (necessarily generic).    Since $\maxim\ul R=\mathbf x\ul R$, \los\ yields that $\seq\maxim w=\seq{\mathbf x}wR$, where $\seq{\mathbf x}w$ is an \sr\ of $\mathbf x$. Since almost all $\seq Rw$ have dimension $d$, it follows that they are regular local rings.

Conversely, suppose almost all $\seq Rw$ are regular. Since the $\seq \br w$-complexity of almost all $\seq Rw$ is at most $c$, for some $c$, we can find a regular system of parameters $\seq{\mathbf x}w$ of $\seq \br w$-complexity at most $c$. By Corollary~\ref{C:gensop}, their ultraproduct $\mathbf x$ is a generic sequence, which generates $\maxim$, by \los\ and faithful flatness. Therefore, $\epsilon\leq d$.   Since always $d\leq\epsilon$ (use Theorem~\ref{T:qdim}), we get that $R$ is \qreg.
\end{proof}

The following is now immediate from the previous result and Theorem~\ref{T:dom}.

\begin{corollary}\label{C:qregdom}
Let $\br$ be an ultra-DVR. If $R$ is an \serene\  \qreg\ local $\br$-affine algebra, then $R$ is a domain and every localization of $R$ with respect to a prime ideal containing $\pi$ is again \qreg.
\end{corollary}

In fact, the \blim\ $R$ of regular local $\seq\br w$-algebras $\seq Rw$ of uniformly bounded $\seq\br w$-complexity is \qreg\ and \serene. Indeed, we already observed that then   $R$ is \serene, and therefore by Theorem~\ref{T:qreg}, \qreg. For a homological characterization of \qreg{ity}, see Corollary~\ref{C:cohreg} below.

\begin{example}
If $R$ denotes the localization of $\pol \br {X,Y}/(X^2+Y^3+\pi)$ at the maximal ideal generated by $X$, $Y$ and $\pi$, then $R$ is \qreg\ (namely $X$ and $Y$ generate the maximal ideal, so $\epsilon =2$, and since $R/\pi R$ has dimension one, $d=2$ as well). Note though that $R/\pi R$ is not regular. 
\end{example}

\subsubsection*{Transfer.}
Let me elaborate on why the results in this section are instances of transfer between positive and mixed \ch. Suppose $\br$ is alternatively realized as the ultraproduct of \DVR{s} $\seq{\tilde \br }w$. Note that this \emph{not} imply that  $\seq \br w$ and $\seq{\tilde \br }w$ are almost all pair-wise isomorphic. In fact, in the next sections, one set of \DVR{s} will be of mixed \ch\ and the other set of prime \ch. Let us put $\seq{\tilde A}w:=\pol{\seq{\tilde \br }w}X$ and let $\ul{\tilde A}$ denote their ultraproduct, so that we have also a canonical embedding $A\to \ul{\tilde A}$. With notation as above, $R$  has also a non-standard hull and \sr{s} with respect to this second set of \DVR{s}; let us denote them by $\ul{\tilde R}$ and $\seq{\tilde R}w$ respectively. Suppose $\seq \br w$ and $\seq{\tilde \br }w$ have pair-wise isomorphic residue fields (as will be the case below). Since the $\seq Rw/\seq\pi w\seq Rw$  are an \sr\ of the $\kappa$-algebra $R/\pi R$ (in the sense of \cite{SchNSTC}) and, mutatis mutandis, so are the $\seq {\tilde R}w/\seq{\tilde\pi} w\seq {\tilde R}w$, where $\seq{\tilde\pi} w$ is a uniformizing parameter of $\seq{\tilde \br }w$, we get from  \cite[3.2.3]{SchNSTC} that almost all $\seq Rw/\seq\pi w\seq Rw$ are isomorphic to $\seq {\tilde R}w/\seq{\tilde\pi} w\seq {\tilde R}w$. Therefore, if we assume that there is no torsion, then $\seq Rw$ and $\seq{\tilde R}w$ have the same dimension, and one set consists of almost all \CM\ local rings \iff\ the other set does (note that this argument does not yet use the above \emph\q\ notions). However, this argument breaks down in the presence of torsion, or, when we want to transfer the regularity property. This can be overcome by using the notions defined in this section, provided we have  a uniform upper bound on the parameter degree.

Suppose, for some $d,e\in\nat$, that almost all $\seq Rw$ have dimension $d$ and parameter degree at most $e$. Note that in view of Corollary~\ref{C:bdpd} this last condition is automatically satisfied if almost all $\seq Rw$ are torsion-free over $\seq\br w$; and that it is implied by the assumption that almost all $\seq Rw$ have uniformly bounded homological multiplicity (see \cite[Corollary 4.6]{SchABCM}).  Applying Proposition~\ref{P:pardeg} twice gives first that $R$ is \serene, with \qdim\ $d$, and then that almost all $\seq {\tilde R}w$ have dimension $d$ and uniformly bounded parameter degree.  Now, Theorems~\ref{T:qCM} and \ref{T:qreg} tell us that almost all $\seq Rw$ are respectively \CM\ or regular, \iff\ almost all $\seq{\tilde R}w$ are.

\section{Big \CM\ Algebras}

In \cite{SchBCM}, ultraproducts of absolute integral closures in \ch\ $p$ were used to define big \CM\ algebras over $\mathbb C$. This same process can be used in the current   mixed \ch\ setting. Recall that for an arbitrary domain $B$, we define its \emph{absolute integral closure} as the integral closure of $B$ in some algebraic closure of its field of fractions and denote it $B^+$. This is uniquely defined up to $B$-algebra isomorphism. 

For each prime number $p$, let $\ulseq \br p{mix}$ be a mixed \ch\ complete \DVR\ with uniformizing parameter $\seq\pi p$ and residue field $\seq \kappa p$ of \ch\ $p$, and let $\br $, $\pi$ and $\kappa$ be their respective ultraproducts. Put $\ulseq \br p{eq}:=\pow{\seq \kappa p}t$, for $t$ a single variable. By the Ax-Kochen-Ershov Theorem, $\br $ is isomorphic to the ultraproduct of the $\ulseq \br p{eq}$. Let $\varpi$ denote  the ideal of infinitesimals of $\br $, that is to say, the intersection of all $\pi^n\br $. As before, we put $A:=\pol \br X$, for a fixed tuple of variables $X$, and let $\uleq A$ and $\ulmix A$ be its respective equi\ch\ and mixed \ch\ non-standard hull, that is to say, the ultraproduct of respectively the $\ulseq Ap{eq}:=\pol{\ulseq \br p{eq}}X$ and the $\ulseq Ap{mix}:=\pol{\ulseq \br p{mix}}X$.

Throughout, $R$ will be a  local $\br$-affine domain with $\ulseq Rp{eq}$ and $\uleq R$ respectively an equi\ch\ \sr\  and the equi\ch\ non-standard hull of $R$ (so that $\uleq R$ is the ultraproduct of the $\ulseq Rp{eq}$).  By Theorem~\ref{T:dom}, almost all $\ulseq Rp{eq}$ are local domains.

\begin{definition}
Define $\BCM R$ as the ultraproduct of the $(\ulseq Rp{eq})^+$.  
\end{definition}

Since $(\ulseq Rp{eq})^+$ is well-defined up to $\ulseq Rp{eq}$-algebra isomorphism, we have that $\BCM R$ is well-defined up to $R$-algebra isomorphism. Moreover, this construction is weakly functorial in the following sense. Let $R\to S$ be an $\br $-algebra \homo\ between local $\br$-affine domains. This induces $\ulseq \br p{eq}$-algebra \homo{s} $\ulseq Rp{eq}\to\ulseq Sp{eq}$ of the corresponding equi\ch\ \sr{s}. These in turn yield \homo{s} $(\ulseq Rp{eq})^+\to (\ulseq Sp{eq})^+$ between the absolute integral closures. Taking ultraproducts, we get an $\br $-algebra \homo\ $\BCM R\to \BCM S$ and a commutative diagram
	\commdiagram [st] R{}S {} {} {\BCM R} {} {\BCM S.}

\begin{theorem}\label{T:BCMmix}
If  $R$  is a local $\br$-affine domain, then any generic sequence in $R$ is $\BCM R$-regular. 
\end{theorem}
\begin{proof}
Let $\uleq R$ and $\ulseq Rp{eq}$ be respectively, the equi\ch\ non-standard hull and an equi\ch\ \sr\ of $R$. Let $\mathbf x$ be a generic sequence, and let $\seq{\mathbf x}p$ be an \sr\ of $\mathbf x$. By Corollary~\ref{C:gensop}, almost all $\seq{\mathbf x}p$ are systems of parameters in $\ulseq Rp{eq}$, whence are $(\ulseq Rp{eq})^+$-regular by \cite{HHbigCM}. By \los, $\mathbf x$ is $\BCM R$-regular.
\end{proof}

\section{Improved New Intersection Theorem}\label{s:INIT}

The remaining sections will establish various asymptotic versions in mixed \ch\ of the Homological Conjectures listed in the abstract. We start with discussing Intersection Theorems. By \cite{Rob87}, we now know that the New Intersection Theorem holds for all Noetherian local rings. However, this is not yet known for the Improved New Intersection Theorem. We need some terminology and notation (all taken from \cite{BH}).

Let $C$ be an arbitrary Noetherian local ring and $\varphi\colon C^a\to C^b$ a linear map between finite free $C$-modules. We will always think of $\varphi$ as  an $(a\times b)$-matrix over $C$. For $r>0$, recall that  the \emph{$r$-th Fitting ideal} of $\varphi$, denoted $I_r(\varphi)$, is the ideal in $C$ generated by all $(r\times r)$ minors of $\varphi$; if $r$ exceeds the size of the matrix, we put $I_r(\varphi):=(0)$. 

With a \emph{finite free complex} over $C$ we mean a complex
	\begin{equation}\label{eq:FFC}
	 0\to C^{a_s}\map{\varphi_s} C^{a_{s-1}}\map{\varphi_{s-1}}\dots \map{\varphi_2} C^{a_1} \map{\varphi_1} C^{a_0}\to 0.\tag{$F_\bullet$}
	\end{equation}
We call $s$ the \emph{length} of the complex, and for each $i$, we define 
	\begin{equation*}
	r_i:= \sum_{j=i}^s (-1)^{j-i} a_j.
	\end{equation*}
We will refer to $r_i$ as the \emph{expected rank} of $\varphi_i$. We will call the residue ring $C/I_{r_i}(\varphi_i)$ the \emph{$i$-th Fitting ring} of $F_\bullet$ and we will denote it $\fitt i(F_\bullet)$.

The $i$-th \emph{homology} of $F_\bullet$ is by definition the quotient module 
	\begin{equation*}
	H_i(F_\bullet):= \op{Ker}(\varphi_i)/\op{Im}(\varphi_{i+1}).
	\end{equation*} 
We call $F_\bullet$ \emph{acyclic}, if all $H_i(F_\bullet)=0$ for $i>0$. In that case, $F_\bullet$ yields a \emph{finite free resolution} of $H_0(F_\bullet)$. 

In case $C$ is $Z$-algebra with $Z$ a Noetherian local ring, we say that $F_\bullet$ has \emph{$Z$-complexity} at most $c$, if its length $s$ is at most $c$, if all $a_i<c$, and if every entry of each $\varphi_i$ has $Z$-complexity at most $c$. Below we will say that an element $\tau$ in a homology module $H_i(F_\bullet)$ has \emph{$Z$-complexity} at most $c$, if it is the image of a tuple in $\op{Ker}(\varphi_i)$  of $Z$-complexity at most $c$ (for more details, see \S\ref{s:tor} below).

\begin{theorem}[Asymptotic Improved New Intersection Theorem]\label{T:INIT}
For each $c$, there exists a bound $\bound{INIT}c$ with the following property.   Let $V$ be a mixed \ch\ \DVR\ and let $(C,\maxim)$ be a $d$-dimensional local $V$-affine domain. Let $F_\bullet$ be a finite free complex over $C$, such that the $i$-th  Fitting ring $\fitt i(F_\bullet)$ has dimension at most $d-i$. Assume $H_0(F_\bullet)$ has  a minimal generator $\tau$   (that is to say, $\tau\notin\maxim H_0(F_\bullet)$), such that $C\tau$ has finite length.

Assume that $c$ simultaneously bounds the $V$-complexity of $C$, $\tau$ and $F_\bullet$, the parameter degree of each Fitting ring $\fitt i(F_\bullet)$, and the length of $C\tau$. If the \ch\ of the residue field of $\br $ is bigger than $\bound{INIT}c$, then $d$ is at most   the  length of the complex $F_\bullet$.
\end{theorem}
\begin{proof}
If $\pi C=0$, then $C$ contains the residue field of $V$ and in that case the Theorem is known (see for instance \cite[Theorem 9.4.1]{BH} or \cite{EG,HoDS}). So we may moreover assume that $C$ is flat over $V$. By faithful flat descent, we may replace $V$ and $C$ by $\complet V$ and $\complet V \tensor_VC$, where $\complet V$ is the completion of $V$. In other words, we only need to prove the result for a torsion-free local domain over a complete \DVR\ of mixed \ch.  Suppose  this last assertion is   false for some $c$, so that for each prime number $p$, we can find a counterexample consisting of the following data:
\begin{itemize}
\item  a mixed \ch\ complete \DVR\ $\ulseq \br p{mix}$ with  uniformizing parameter $\seq\pi p$, whose residue field has \ch\ $p$;
\item  a  local  $\ulseq \br p{mix}$-affine domain $\ulseq Rp{mix}$ of $\ulseq \br p{mix}$-complexity at most $c$;
\item a finite free complex
	\begin{multline}\label{eq:FFCp}
	 0\to (\ulseq Rp{mix})^{a_s}\map{\seq{\varphi_{s,}}p} (\ulseq Rp{mix})^{a_{s-1}}\map{\seq{\varphi_{s-1,}}p}\dots \\
\map{\seq{\varphi_{2,}}p} (\ulseq Rp{mix})^{a_1} \map{\seq{\varphi_{1,}}p} (\ulseq Rp{mix})^{a_0}\to 0.\tag{$F_{p\bullet}^{\text{mix}}$}
	\end{multline}
of length $s$ and of $\ulseq \br p{mix}$-complexity at most $c$, such that the $i$-th Fitting ring $\fitt i(\ulseq F{p\bullet}{mix})$ has dimension at most $d-i$ and parameter degree at most $c$;
\item a minimal generator $\seq \tau p$ of $H_0(\ulseq F{p\bullet}{mix})$ of $\ulseq \br p{mix}$-complexity at most $c$, generating a module of length at most $c$, 
\end{itemize}
but such that $s$ is strictly less than the dimension of $\ulseq Rp{mix}$. Note that, without loss of generality, we may assume that the dimension of $\ulseq Rp{mix}$ and that all the ranks of $F_{p\bullet}^{\text{mix}}$  are independent from $p$, since there are only finitely many possibilities, so that precisely one such possibility almost always holds. In particular, the expected ranks do not depend on $p$.

 Let $\br $ and $\pi$ be the respective ultraproduct of the $\ulseq \br p{mix}$ and the $\seq \pi p$. Let $R$ and $\ulmix R$ be the respective \blim\ and ultraproduct of the $\ulseq Rp{mix}$. It follows from Theorem~\ref{T:dom}, that $R$ is a local $\br$-affine domain, and from Theorem~\ref{T:ff}, that $R\to \ulmix R$ is faithfully flat. Let $d$ be the \qdim\ of $R$, so that almost all $\ulseq Rp{mix}$ have dimension $d$ by Theorem~\ref{T:qdim}. Let $\varphi_i$ be the ultraproduct of the $\seq{\varphi_{i,}}p$. It follows from Lemma~\ref{L:deg} that each $\varphi_i$ is already defined over $R$. Hence by \los, we get that
	\begin{equation}\label{eq:FFCR}
	 0\to R^{a_s}\map{\varphi_s} R^{a_{s-1}}\map{\varphi_{s-1}}\dots \map{\varphi_2} R^{a_1} \map{\varphi_1} R^{a_0}\to 0\tag{$F_\bullet$}
	\end{equation}
is a finite free complex. Let $M$ denote its zero-th homology. Fix some $i$. By \los, $I_{r_i}(\seq{\varphi_{i,}}p)$  is an $\ulseq \br p{mix}$-\sr\ of $I_{r_i}(\varphi_i)$. By the uniform boundedness of the parameter degrees, $\fitt i(F_\bullet)$ is \serene\ by Proposition~\ref{P:pardeg}. If $d_i$ is the \qdim\ of $\fitt i(F_\bullet)$, then $d-d_i$ is equal to the height of almost all $I_{r_i}(\seq{\varphi_{i,}}p)$ and to the  \qht\ of $I_{r_i}(\varphi_i)$, by Theorem~\ref{T:qht}. In particular, by assumption, $i\leq d-d_i$, and therefore, by definition of \qht, we can find a generic sequence  $\mathbf x_i$ in $R$ whose first $i$ entries belong to $I_{r_i}(\varphi_i)$. 

Let $B:=\BCM R$. Since $\mathbf x_i$ is $B$-regular by Theorem~\ref{T:BCMmix}, the grade of  $I_{r_i}(\varphi_i)B$ is at least $i$. Since this holds for all $i$, the Buchsbaum-Eisenbud-Northcott Acyclicity Theorem  (\cite[Theorem 9.1.6]{BH}) proves that $F_\bullet\tensor_RB$ is acyclic. Since $B$ has depth at least  $d$, it follows from \cite[Theorem 9.1.2]{BH} that  the zero-th homology of $F_\bullet\tensor_RB$, that is to say, $M\tensor_RB$, has depth at least $d-s$. 

Let $\tau $ be the ultraproduct of the $\seq \tau p$. Note that each $\seq \tau p$ is by assumption the image of a tuple in $(\ulseq Rp{mix})^{a_0}$ of $\ulseq \br p{mix}$-complexity at most $c$, so that $\tau $ is already defined over $R$ by Lemma~\ref{L:deg}. By \los, $\tau $ is a minimal generator of 
	\begin{equation*}
	H_0(F_\bullet\tensor\ulmix R)=M\tensor\ulmix R,
	\end{equation*}
and by \cite[Proposition 1.1]{SchEC} or \cite[Proposition 9.1]{JL89}, the length of $\ulmix R\tau $ is at most $c$. By faithful flatness, $\tau \in M-\maxim M$, where $\maxim$ is the maximal ideal of $R$, and $R\tau $ has length at most $c$. In particular, the image of $\tau \tensor 1$ in $M/\maxim M\tensor B/\maxim B$ is non-zero, and therefore $\tau \tensor 1$ itself is a non-zero element of $M\tensor B$. Since   $\maxim^c$ annihilates $\tau \tensor 1$, we get that $M\tensor B$ has depth  zero. Together with the conclusion from the previous paragraph, we get that $d\leq s$, contradiction. 
\end{proof}

\section{Monomial and Direct Summand Conjectures}

We keep notation as in the previous section, so that in particular $\br $ will denote the ultraproduct of mixed \ch\ complete \DVR{s} $\ulseq \br p{mix}$. In order to formulate a non-standard version of the  Monomial Conjecture, we need some terminology. Let $\ul\nat$ be the ultrapower of $\nat$. Let $\seq Cw$ be rings, $X:=\rij Xd$ variables and $\ul A$ the ultraproduct of the $\pol{\seq Cw}X$. Although each $\pol{\seq Cw}X$ is $\nat$-graded, it is not true that $\ul A$ is $\ul\nat$-graded, since we might have infinite sums of monomials in $\ul A$. Nonetheless, for each $\ul\nu\in\ul\nat$, the element $X^{\ul\nu}$ is well-defined, namely, if $\ul\nu$ is the ultraproduct of elements $\seq\nu w\in \nat$, then
	\begin{equation*}
	X^{\ul\nu}:= \up w X^{\seq\nu p}.
	\end{equation*}
In particular, if $\ul B$ is an arbitrary ultraproduct of rings $\seq Bw$ and if $\mathbf x$ is a $d$-tuple in $\ul B$, then $\mathbf x^{\ul\nu}$ is a well-defined element of $\ul B$.

With a \emph{cone} $H$ in a semi-group $\Gamma$ (e.g.,  $\Gamma=\nat^d$ or $\Gamma=\ul\nat^d$), we mean a subset $H$ of $\Gamma$ such that $\nu+\gamma\in H$, for every $\nu\in H$ and every $\gamma\in\Gamma$. A cone $H$ is \emph{finitely generated}, if there exist $\nu_1,\dots,\nu_s\in H$, called \emph{generators} of the cone, such that
	\begin{equation*}
	H= \bigcup_i\nu_i+\Gamma.
	\end{equation*}
If $H$ is a cone in $\nat^d$, we denote by $J_H$ the monomial ideal in $\pol\zet Y$, generated by all $Y^\nu$ with $\nu\in H$, where $Y$ is a $d$-tuple of variables. If $H$ is generated by $\nu_1,\dots,\nu_s$, then $J_H$ is generated by $X^{\nu_1},\dots, X^{\nu_s}$. Conversely, if $J$ is a monomial ideal in $\pol\zet Y$, then the collection of all $\nu$ for which $Y^\nu\in J$, is a cone in $\nat^d$. Since $\pol\zet Y$ is Noetherian,   every cone in $\nat^d$ is finitely generated. This is no longer true for a cone in $\ul\nat^d$.

Let $B$ be an arbitrary ring. We will use the following well-known fact about regular sequences. If $\mathbf x$ is a $B$-regular sequence (in fact, it suffices that $\mathbf x$ is quasi-regular), $H$ a cone in $\nat^d$ and $\nu\notin H$, then $\mathbf x^\nu$ does not lie in the ideal $J_H(\mathbf x)$ generated by all $\mathbf x^\theta$ with $\theta\in H$.

\begin{corollary}\label{C:nsMC}
Let  $R$  be a local $\br$-affine domain with equi\ch\ non-standard hull $\uleq R$. Let  $\mathbf x$  be a generic sequence in $R$, let $H$ be a cone in $\ul\nat^d$ and let $\ul\nu\in\ul\nat^d$. If   $\ul\nu\notin  H$, then
	\begin{equation}\label{eq:nsmc}
	\mathbf x^{\ul\nu} \notin J_H(\mathbf x):=(\mathbf x^\mu\mid\mu\in H)\uleq R.
	\end{equation}
\end{corollary}
\begin{proof}
Suppose \eqref{eq:nsmc} is false for some choice of cone $H$ of $\ul\nat^d$ and some $\ul{\nu_0}\notin H$. In other words, we can find $\ul{f_i}$ in $\uleq R$ and tuples $\ul{\nu_i}$ in $H$, such that
	\begin{equation}\label{eq:nsmcc}
	\mathbf x^{\ul{\nu_0}} = \ul{f_1} \mathbf x^{\ul{\nu_1}}+ \dots + \ul{f_s} \mathbf x^{\ul{\nu_s}}.
	\end{equation}
In order to derive a contradiction, we will argue that such a relation~\eqref{eq:nsmcc} cannot hold in $\BCM R$. Indeed, suppose it does. Let   $\ulseq Rp{eq}$ be an equi\ch\ \sr\ of $R$, so that $\BCM R$ is the ultraproduct of the $(\ulseq Rp{eq})^+$. Choose tuples $\seq {\nu_i}p\in\nat$, elements $\seq{f_i}p\in(\ulseq Rp{eq})^+$ and tuples $\seq{\mathbf x}p$ in $\ulseq Rp{eq}$ whose respective ultraproducts are $\ul{\nu_i}$, $\ul{f_i}$ and $\mathbf x$. By \los, we get that
	\begin{equation}\label{eq:nsmcp}
	\seq{\mathbf x}p^{\seq{\nu_0}p} = \seq{f_1}p \seq{\mathbf x}p^{\seq{\nu_1}p}+ \dots + \seq{f_s}p\seq{\mathbf x}p^{\seq{\nu_s}p}
	\end{equation}
in $(\ulseq Rp{eq})^+$, for almost all $p$. \los\ also yields that $\seq{\nu_0}p$ does not lie in the cone of $\nat^d$ generated by $\seq{\nu_1}p,\dots,\seq{\nu_s}p$, for almost all $p$. However,  $\mathbf x$ is $\BCM R$-regular by Theorem~\ref{T:BCMmix}, whence, almost all $\seq{\mathbf x}p$ are $(\ulseq Rp{eq})^+$-regular by \los. It follows that \eqref{eq:nsmcp} cannot hold for those $p$.
\end{proof}

\begin{theorem}[Asymptotic Monomial Conjecture]\label{T:MC}
For each pair $(c,s)$, there exists a bound $\bound {MC}{c,s}$ with the following property. Let $Y$ be a $d$-tuple of variables, $J$ a monomial ideal  in $\pol\zet Y$ generated by $s$ monomials and $Y^\nu$ a monomial not belonging to $J$. Let $V$ be a mixed \ch\ \DVR\  and let $C$ be a $d$-dimensional local $V$-affine domain. Let $\mathbf y$ be a system of parameters  in $C$.  

If $C$ and $\mathbf y$ have $V$-complexity at most $c$ and $\pi^c\in\mathbf yC$, and if the \ch\ of the residue field of $V$ is bigger than  $\bound {MC}{c,s}$, then 
	\begin{equation*}
	\mathbf y^\nu\notin J(\mathbf y)C
	\end{equation*}
where $J(\mathbf y)C$ is the ideal in $C$ obtained from $J$ by the substitution $Y\mapsto \mathbf y$.
\end{theorem}
\begin{proof}
Note that since  $C$ has $V$-complexity at most $c$, its dimension $d$ is at most $c$. By faithful flat descent, we may replace $V$ and $C$ by $\complet V$ and $\complet V\tensor_VC$, where $\complet V$ is the completion of $V$. In other words, we only need to prove the result for complete \DVR{s} of mixed \ch. Suppose this is false for some $(c,s)$, so that we can find for each prime number $p$
\begin{itemize}
\item a mixed \ch\ complete \DVR\ $\ulseq \br p{mix}$ with  uniformizing parameter $\seq\pi p$, whose residue field has \ch\ $p$,
\item a   local $\ulseq \br p{mix}$-affine domain $\ulseq Rp{mix}$ of $\ulseq \br p{mix}$-complexity at most $c$,
\item  tuples   $\seq{\nu_0}p,\dots,\seq{\nu_s}p$ with $\seq{\nu_0}p$ not in the cone generated by the remaining tuples,
\item  a system of parameters $\seq{\mathbf y}p$ of $\ulseq \br p{mix}$-complexity at most $c$ generating an ideal containing $(\seq\pi p)^c$, 
\end{itemize}
such that
	\begin{equation*}
	\seq{\mathbf y}p^{\seq{\nu_0}p} \in (\seq{\mathbf y}p^{\seq{\nu_1}p},\dots, \seq{\mathbf y}p^{\seq{\nu_s}p})\ulseq Rp{mix}.
	\end{equation*}
Let $\br $ be the ultraproduct of the $\ulseq \br p{mix}$ and let $R$ and $\ulmix R$ be the respective \blim\ and ultraproduct of the $\ulseq Rp{mix}$. Since $R$ is then a domain, it  is \serene. Let $\mathbf y$ and $\ul{\nu_i}$ be the respective ultraproducts of $\seq{\mathbf y}p$ and $\seq{\nu_i}p$.  The sequence $\mathbf y$ is defined over $R$, by Lemma~\ref{L:deg}, and is generic in $R$, by Corollary~\ref{C:gensop}. By \los\ and Theorem~\ref{T:ff}, we get that
	\begin{equation*}
		 \mathbf y^{\ul{\nu_0}} \in ( \mathbf y^{\ul{\nu_1}},\dots, \mathbf y^{\ul{\nu_s}})R.
	\end{equation*} 
However, this contradicts Corollary~\ref{C:nsMC} for $H$ the cone of $\ul\nat^d$ generated by $\ul{\nu_1},\dots,\ul{\nu_s}$.
\end{proof}

\begin{remark}
Using some results from \cite{SchUBS}, we can remove the restriction on $C$ to be a domain. Namely, by the usual argument, we reduce to the domain case by killing a minimal prime $\pr$ of $C$  of maximal dimension (that is to say, so that $\op{dim}C=\op{dim}C/\pr$). However, in order to apply the theorem to the domain $C/\pr$, we must be guaranteed that its   $V$-complexity is at most $c'$, for some $c'$ only depending on $c$. Such a bound does indeed exists by \cite[Theorems 9.2 and 9.12]{SchUBS}.
\end{remark}

\begin{theorem}[Asymptotic Direct Summand Conjecture]
For each $c$, we can find a bound $\bound{DS}c$ with the following property. Let $V$ be a mixed \ch\ \DVR\  and let $C\to D$ be a finite injective local \homo\ between local $V$-affine algebras, with $C$ regular.

If $C\to D$ has $V$-complexity at most $c$ and if the \ch\ of the residue field of $V$ is bigger than the bound $\bound{DS}c$, then $C$ is a direct summand of $D$ (as a $C$-module).
\end{theorem}
\begin{proof}
If $\pi C=0$, we are in the equi\ch\ case and the result is well-known. So we may assume that $V\subset C$. We leave it to the reader to make the reduction to the case that $V$ is complete and $D$ is torsion-free over $V$. Choose a regular system of parameters $\mathbf x:=\rij xd$ of $C$ of $V$-complexity at most $c$. Since $C\to D$ is finite, $\mathbf x$ is a system of parameters in $D$ with  $\pi\in\mathbf xD$. Since $C\to D$ has $V$-complexity at most $c$, the image of $\mathbf x$ in $D$ has $V$-complexity at most $c^2$. Hence if the \ch\ of the residue field is bigger than $\bound{MC}{c^2,d}$, then by Theorem~\ref{T:MC}
	\begin{equation*}
	(x_1x_2\cdots x_d)^t \notin (x_1^{t+1},\dots, x_d^{t+1})D,
	\end{equation*}
for any $t$. By \cite[Lemma 9.2.2]{BH}, this implies that $C$ is a direct summand of $D$.
\end{proof}

\section{Pure subrings of regular rings}

We keep notation as in the previous section, so that in particular $\br $ will denote the ultraproduct of mixed \ch\ complete \DVR{s} $\ulseq \br p{mix}$. Our goal is to show an asymptotic version of the Hochster-Roberts Theorem in \cite{HR}. Recall that a ring \homo\ $C\to D$ is called \emph{cyclically pure} if every ideal $I$ in $C$ is extended from $D$, that is to say, if $I=ID\cap C$.

\begin{theorem}\label{T:BCMreg}
If $R$ is a \qreg\ \serene\ local $\br$-affine algebra, then $R\to\BCM R$ is faithfully flat.
\end{theorem}
\begin{proof}
Let $L$ be a linear form in a finite number of variables $Y$ with coefficients in $R$ and let $\mathbf b$ be a solution in $B:=\BCM R$ of $L=0$. Let $\ulseq Rp{eq}$, $\ulseq Lp{eq}$ and $\ulseq{\mathbf b}p{eq}$ be equi\ch\ \sr{s} of $R$, $L$ and $\mathbf b$ respectively. By \los, $\ulseq{\mathbf b}p{eq}$ is a solution in $(\ulseq Rp{eq})^+$ of the linear equation $\ulseq Lp{eq}=0$. By \cite[Corollary 4.27]{Asch}, we can find tuples $\ulseq{\mathbf a_1}p{eq}, \dots, \ulseq{\mathbf a_s}p{eq}$ over $\ulseq {R}p{eq}$ generating the module of solutions of $\ulseq {L}p{eq}=0$, all of $\ulseq {\br}p{eq}$-complexity at most $c$, for some $c$ independent from $p$ and $s$. Let $\mathbf a_1,\dots,\mathbf a_s$ be the respective ultraproducts, which are then defined over $R$ by Lemma~\ref{L:deg}. By \los, $L(\mathbf a_i)=0$, for each $i$. On the other hand, almost all $\ulseq {R}p{eq}$ are regular, by Theorem~\ref{T:qreg}. Therefore, $\ulseq {R}p{eq}\to (\ulseq {R}p{eq})^+$ is flat by \cite[Theorem 9.1]{HuTC}. Hence we can write $\ulseq {\mathbf b}p{eq}$ as a linear combination over $(\ulseq {R}p{eq})^+$ of the $\ulseq{\mathbf a_i}p{eq}$. By \los, $\mathbf b$ is a $B$-linear combination of the solutions $\mathbf a_i$, showing that $R\to B$ is flat whence faithfully flat.
\end{proof}

\begin{proposition}\label{P:psHR}
Let $R\to S$ be an injective \homo\ of local \serene\ $\br$-affine algebras. If $R/\pi R\to S/\pi S$ is cyclically pure  and $S$ is a  \qreg\ local ring, then $R$ is \qCM.
\end{proposition}
\begin{proof}
Since $S$ is a domain by Corollary~\ref{C:qregdom}, so is $R$ by cyclic purity. If $\pi R=0$, we are in an equi\ch\ Noetherian situation and the statement becomes the  Hochster-Roberts Theorem \cite{HR}.   Therefore, we may assume $\pi$ is $R$-regular, so that we can choose a generic sequence $\mathbf x:=\rij xd$ in $R$ with  $x_1=\pi$. For each $n$, let $I_n:=\rij xnR$. Suppose $rx_{n+1}\in I_n$, for some $r\in R$. By Theorem~\ref{T:BCMmix}, we have that $\mathbf x$ is a $\BCM R$-regular sequence. Therefore, $r\in I_n\BCM R$. Since the \homo\ $R\to S$ induces a \homo\  $\BCM R\to \BCM S$, we get that $r\in I_n\BCM S$.  By Theorem~\ref{T:BCMreg}, we have $I_n\BCM S\cap S=I_nS$, so that $r\in I_nS$. Using finally that $R/\pi R\to S/\pi S$ is cyclically pure, we get $r\in I_n$. This shows that $\mathbf x$ is $R$-regular, so that $R$ has depth at least $d$ and hence is \qCM. 
\end{proof}

\begin{theorem}[Asymptotic Hochster-Roberts Theorem]
For each $c$, we can find a bound $\bound{HR}c$ with the following property. Let $V$ be a mixed \ch\ \DVR\ and let $C\to D$ be a cyclically pure local  \homo\  of local $V$-algebras with $D$ regular. 

If $C\to D$ has complexity at most $c$, then $C$ is \CM, provided the \ch\ of the residue field of $V$ is at least $\bound{HR}c$.
\end{theorem}
\begin{proof}
As before, we may reduce to the case that $V$ is complete and that $V\subset C$. Suppose this assertion is then false for some $c$, so that we can find for each prime number $p$, a mixed \ch\ complete \DVR\ $\ulseq \br p{mix}$ with residue field of \ch\ $p$  and a cyclically pure \homo\ $\ulseq Rp{mix}\to\ulseq Sp{mix}$  of $\ulseq \br p{mix}$-complexity at most $c$ between  local $\ulseq \br p{mix}$-algebras, such that $\ulseq Sp{mix}$ is regular but $\ulseq Rp{mix}$ is not \CM. Let  $R\to S$ and  $\ulmix R\to \ulmix S$ be respectively the \blim\ and the ultraproduct of   the $\ulseq Rp{mix}\to \ulseq Sp{mix}$. Theorem~\ref{T:qCM} implies that $R$ is not \qCM\ and Theorem~\ref{T:qreg}, that $S$ is \qreg. I claim that $R/\pi R\to S/\pi S$  is cyclically pure. Assuming this claim, we get from Proposition~\ref{P:psHR} that $R$ is \qCM, contradiction.

To prove the claim, let $I$ be an arbitrary ideal in $R$ containing $\pi$. Let  $r\in I S\cap  R$, so that  we need to show that $r\in I$. Note that $I$ is   finitely generated, as $R/\pi R$ is Noetherian. Let $\ulseq Ip{mix}$ and $\ulseq rp{mix}$ be mixed \ch\ \sr{s} in $\ulseq Rp{mix}$ of $I$ and $r$ respectively. By \los, almost all $\ulseq rp{mix}$ lie in $\ulseq Ip{mix}\ulseq Sp{mix}\cap\ulseq Rp{mix}$, whence in $\ulseq Ip{mix}$ by cyclical purity.  Hence by \los\ $r\in I\ulmix R$, so that  $r\in I$ by faithful flatness, as we needed to prove.
\end{proof}

\section{Asymptotic vanishing for maps of Tor}\label{s:tor}

\begin{proposition}\label{P:int}
If $R\to S$ is an  integral extension of local $\br$-affine domains, then $\BCM R=\BCM S$.
\end{proposition}
\begin{proof}
Since any integral extension is a direct limit of finite extensions, we may assume that $R\to S$ is finite. Choose an equi\ch\ \sr\ $\ulseq Rp{eq}\to \ulseq Sp{eq}$ of $R\to S$. By Theorem~\ref{T:dom} and \los, almost all $\ulseq Rp{eq}$ and $\ulseq Sp{eq}$ are domains and the extension $\ulseq Rp{eq}\to\ulseq Sp{eq}$ is finite. Therefore,   $(\ulseq Rp{eq})^+=(\ulseq Sp{eq})^+$, so that in the ultraproduct, we get $\BCM R=\BCM S$.
\end{proof}

\begin{theorem}\label{T:qtor}
Let $R\to S\to T$ be local $\br$-algebra \homo{s} between local $\br$-affine domains. Assume that $R$ and $T$ are \qreg\ and that $R\to S$ is integral and injective. For every  $R$-module $M$, the induced map $\tor RiSM\to\tor RiTM$ is zero, for all $i\geq 1$.
\end{theorem}
\begin{proof}
Since $R\to S$ is integral, we have that $\BCM R=\BCM S$ by Proposition~\ref{P:int}. Therefore, $\tor Ri{\BCM S}M=0$, for all $i\geq 1$, by Theorem~\ref{T:BCMreg}. By weak functoriality, we have, for each $i\geq 1$, a commutative diagram
	\commdiagram {\tor RiSM} {} {\tor RiTM} {} {} {0=\tor Ri{\BCM S}M} {} {\tor Ri{\BCM T}M.}
In particular, the composite map in this diagram is zero, so that  the statement follows once we have shown that the last vertical map is injective. However, this is clear, since $T\to \BCM T$ is faithfully flat by Theorem~\ref{T:BCMreg}. 
\end{proof}

To make use of this theorem, we need to incorporate modules in our present  setup. I will not provide full details, since many results are completely analogous to the case where we work over a field, and this has been treated in detail in \cite{SchBC}. Of course, we do not have the full equivalent of Theorem~\ref{T:SvdD} to our disposal, but for most purposes, the flatness result in Theorem~\ref{T:ff} suffices.

Let $C$ be an arbitrary Noetherian local ring and $M$ a finitely generated module over $C$. We say that a finite free complex $F_\bullet$ is a \emph{finite free resolution of $M$ up to level $n$}, if   $H_0(F_\bullet)=M$ and all $H_j(F_\bullet)=0$, for $j=\range 1n$. Hence, if $n$ is strictly larger than the length of $F_\bullet$, then this just means that $F_\bullet$ is a finite free resolution of $M$ (compare with the terminology introduced in the beginning of \S\ref{s:INIT}). 

Suppose moreover that $Z$ is a Noetherian local ring and $C$ is a local $Z$-affine algebra. We say that $M$ has \emph{$Z$-complexity} at most $c$, if $C$ has $Z$-complexity at most $c$ and if $M$ can be realized as the cokernel of a matrix of $Z$-complexity at most $c$ (meaning that its size is at most $c$ and all its entries have $Z$-complexity at most $c$).

\begin{proposition}\label{P:res}
For each pair $(c,n)$, there exist   bounds $\bound{RES}{c,n}$ and $\bound{HOM} c$ with the following property. Let $V$ be a mixed \ch\ \DVR\ and let $C$ be a local $Z$-affine algebra of $V$-complexity at most $c$.  
\begin{itemize}
\item Any   finitely generated $C$-module of $V$-complexity at most $c$, admits a (minimal) finite free resolution up to level $n$ of $V$-complexity at most $\bound{RES}{c,n}$.
\item Any finite free complex over $C$ of $V$-complexity at most $c$, has homology modules of $Z$-complexity at most $\bound{HOM}c$.
\end{itemize}
\end{proposition}
\begin{proof}
The first assertion follows by induction from the already quoted \cite[Corollary 4.27]{Asch} on bounds of syzygies (compare with the proof of \cite[Theorem 4.3]{SchBC}). It is also clear that we may take this resolution to be minimal (=every tuple in one of the kernels has its entries in the maximal ideal), if we choose to do so. The second assertion is derived from the flatness of the non-standard hull in exactly the same manner as the corresponding result for fields was obtained in \cite[Lemma 4.2 and Theorem 4.3]{SchBC}.
\end{proof}

Recall that the \emph{weak global dimension}  of a ring $C$ is by definition the supremum (possibly infinite) of the weak homological dimensions (=flat dimensions) of all $C$-modules, that is to say, the supremum of all $n$ for which $\tor Cn\cdot\cdot$ is not identically zero.

\begin{corollary}\label{C:wgdim}
A   \qreg\ local $\br$-affine domain  has finite weak global dimension.
\end{corollary}
\begin{proof}
Let $R$ be a \qreg\ local $\br$-affine domain. Given an arbitrary $R$-module $M$, we have to show that $M$ has finite flat dimension, that is to say, admits a finite flat resolution. Assume first that $M$ is finitely presented. Hence we can realize $M$ as the cokernel of some  matrix $\Gamma$. Let $\ul R$ be the non-standard hull of $R$ and let $\seq Rw$ and $\seq \Gamma w$ be \sr{s} of $R$ and $\Gamma$ respectively. Let $\seq Mw$ be the cokernel of $\seq\Gamma w$.  Let $d$ be the  \qdim\ of $R$. By Proposition~\ref{P:res}, we can find a finite free resolution $F_{w\bullet}$ up to level $d$ of each $\seq Mw$, of $\seq\br w$-complexity at most $c$, for some $c$ depending only on $\Gamma$, whence independent from $w$. Since almost all $\seq Rw$ are regular by Theorem~\ref{T:qreg} and have dimension $d$ by Theorem~\ref{T:qdim}, each $\seq Mw$ has projective dimension at most $d$, so that we can even assume that $F_{w\bullet}$ is a finite free resolution of $\seq Mw$. Let $F_\bullet$ be the \blim\ of the $F_{w\bullet}$ (that is to say, the finite free complex over $R$ given by the \blim\ of the matrices in $F_{w\bullet}$). By \los, $F_\bullet\tensor_R\ul R$ is a free resolution of $M\tensor_R\ul R$, and therefore by faithful flat descent, $F_\bullet$ is a free resolution of $M$, proving that $M$ has projective dimension at most $d$.

Assume now that $M$ is arbitrary.  By what we just proved, we have for every finitely generated ideal $I$ of $R$  that $\tor R{d+1}M{R/I}$ vanishes. Hence, if $H$ is a $d$-th syzygy of $M$, then $\tor R1H{R/I}=0$. Since this holds for every finitely generated ideal of $R$, we get from \cite[Theorem 7.7]{Mats} that $H$ is flat over $R$. Hence $M$ has finite flat dimension (at most $d$). 
\end{proof}

By \cite{Jen70}, any flat $R$-module has projective dimension less than the finitistic global dimension of $R$ (the supremum of all projective dimensions of modules of finite projective dimension). Therefore, if, moreover,  the finitistic global dimension of $R$ is finite, then so is its global dimension.  For a Noetherian local ring, its global dimension is finite \iff\ its residue field has finite projective dimension (\iff\ it is regular). The following is the \q\ analogue of this.

\begin{corollary}\label{C:cohreg}
A    local $\br$-affine domain  is \qreg\ \iff\ it is a coherent regular ring in the sense of \cite{Bert}  \iff\ its residue field has finite projective dimension.
\end{corollary}
\begin{proof}
In \cite{Bert} or \cite[\S5]{Glaz}, a local ring $R$ is called a \emph{coherent regular ring}, if every finitely generated ideal of $R$ has finite projective dimension. If $R$ is \qreg\ local $\br$-affine domain, then this property was established in the course of the proof of Corollary~\ref{C:wgdim}. Conversely, suppose $R$ is a local $\br$-affine domain in which every finitely generated ideal has finite projective dimension. In particular, its residue field $k$ admits a finite projective resolution, say of length $n$. Let $\seq Rw$ and $\seq kw$ be \sr{s} of $R$ and $k$ respectively.  Since the $\seq kw$ have uniformly bounded $\seq\br w$-complexity, Proposition~\ref{P:res} allows us to take a minimal finite free resolution $F_{w\bullet}$ of $\seq kw$ up to level $n$, with the property that each $F_{w\bullet}$ has  $\seq\br w$-complexity at most $c$, for some $c$ independent from $w$. Let $F_\bullet$ be the \blim\ of these resolutions. By \los, $F_\bullet$ is a minimal finite free resolution of $k$ up to level $n$. Since $F_\bullet$ is minimal and since $k$ has by assumption projective dimension $n$, it follows that the final morphism (that is to say, the left most arrow)  in $F_\bullet$ is injective. By \los, so are almost all final morphisms in $F_{w\bullet}$, showing that almost all $\seq kw$ have finite projective dimension. By Serre's characterization of regular local rings, we conclude that almost all $\seq Rw$ are regular. Theorem~\ref{T:qreg} then yields that $R$ is \qreg, as we wanted to show.
\end{proof}

Closer inspection of the above argument shows that the residue field of a \qreg\   local $\br$-affine domain $R$ has projective dimension equal to the \qdim\ of $R$. In particular,  the weak global dimension of $R$ is equal to its \qdim. 

\begin{theorem}[Asymptotic Vanishing for Maps of Tors]
For each $c$, we can find a bound $\bound{VT}{c}$ with the following property. Let $V$ be a mixed \ch\ \DVR,  let $C\to D\to E$ be   local  $V$-algebra \homo{s}  of local domains and let $M$ be a finitely generated $R$-module. Assume $C$ and $E$ are regular and $C\to D$ is integral and injective.

If $M$ and all \homo{s} have $V$-complexity at most $c$, then the natural map $\tor CnDM\to\tor CnEM$ is zero, for all $n\geq 1$, provided the \ch\ of the residue field of $V$ is at least $\bound{VT}{c}$.
\end{theorem}
\begin{proof}
Note that $C$ has dimension at most $c$ and therefore $\tor Cn\cdot\cdot$ vanishes identically for all $n>c$ and the assertion trivially holds for these values of $n$. If $\pi C=0$, we are in the equi\ch\ case and the result is known in that case (\cite[Theorem 9.7]{HuTC}). Hence we may assume that all rings are torsion-free over $V$. Moreover, without loss of generality, we may assume that $V$ is complete. Suppose even in this restricted setting, there is no such bound for $c$ and some $1\leq n<c$. Hence, for each prime number $p$, we can find a counterexample consisting of the following data:
\begin{itemize}
\item  a mixed \ch\ complete \DVR\ $\ulseq \br p{mix}$   of residual \ch\ $p$;
\item   local \homo{s} $\ulseq Rp{mix}\to \ulseq Sp{mix}\to \ulseq Tp{mix}$ of torsion-free local domains  of $\ulseq \br p{mix}$-complexity at most $c$, with $\ulseq Rp{mix}$ and $\ulseq Tp{mix}$ regular and $\ulseq Rp{mix}\to \ulseq Sp{mix}$ integral;
\item a finitely generated $\ulseq Rp{mix}$-module $\ulseq Mp{mix}$ of $\ulseq \br p{mix}$-complexity at most $c$; 
\end{itemize}
such that 
	\begin{equation*}
	\tor {\ulseq Rp{mix}} n {\ulseq Sp{mix}}{\ulseq Mp{mix}} \to \tor {\ulseq Rp{mix}} n {\ulseq Tp{mix}}{\ulseq Mp{mix}}
	\end{equation*}
is non-zero.

Let $\br$ be the ultraproduct of the $\ulseq \br p{mix}$ and let $M$ be the \blim\ of the $\ulseq Mp{mix}$ (that is to say, $M$ is the cokernel of the \blim\ of matrices whose cokernel is $\ulseq Mp{mix}$). Let $R\to S\to T$ and $\ulmix R\to\ulmix S\to \ulmix T$ be the respective \blim\ and mixed \ch\ ultraproduct of the \homo{s} $\ulseq Rp{mix}\to \ulseq Sp{mix}\to \ulseq Tp{mix}$. It follows from Corollary~\ref{C:qregdom} and Theorems~\ref{T:dom} and \ref{T:qreg}, that $R$, $S$ and $T$ are local $\br$-affine domains with $R$ and $T$ \qreg. By \los, $\ul R\to \ul S$ is integral, whence so is $R\to S$ by faithful flat descent. By Theorem~\ref{T:qtor}, the natural \homo\ $\tor RnSM\to \tor RnTM$ is therefore zero.

By Proposition~\ref{P:res}, we can find  a finite free resolution  $F_{p\bullet}^{\text{mix}}$ of $\ulseq Mp{mix}$ up to level $n$, of $\ulseq\br p{mix}$ complexity at most $c'$, for some $c'$ only depending on $c$ and $n$. By definition of Tor, we have isomorphisms
	\begin{align*}
	\tor {\ulseq Rp{mix}} n {\ulseq Sp{mix}}{\ulseq Mp{mix}} &\iso H_n(F_{p\bullet}^{\text{mix}}\tensor_{\ulseq Rp{mix}} \ulseq Sp{mix})\\
	\tor {\ulseq Rp{mix}} n {\ulseq Tp{mix}}{\ulseq Mp{mix}} &\iso H_n(F_{p\bullet}^{\text{mix}}\tensor_{\ulseq Rp{mix}} \ulseq Tp{mix})
	\end{align*}
In particular, by Proposition~\ref{P:res}, both modules have $\ulseq \br p{mix}$-complexity at most $c''$, for some $c''$ only depending on $c'$, whence only on $(c,n)$. Let $H_S$ and $H_T$ be their respective \blim, so that by \los\ and our assumptions, $H_S\to H_T$ is non-zero. Let $F_\bullet$ be the \blim\ of the $F_{p\bullet}^{\text{mix}}$. By \los\ and faithful flatness,   $H_S$ and $H_T$ are isomorphic respectively with $H_n(F_\bullet\tensor_RS)$ and $H_n(F_\bullet\tensor_RT)$. However, since $F_\bullet$ is a finite free resolution of $M$ up to level $n$ by another application of \los\ and faithful flatness, these two modules are  isomorphic respectively with $\tor RnSM$  and $\tor RnTM$. Hence the natural map between these two modules is non-zero, contradiction.
\end{proof}

\providecommand{\bysame}{\leavevmode\hbox to3em{\hrulefill}\thinspace}

\end{document}